\documentclass[11pt]{amsart}
\usepackage{amsmath,amssymb,amsthm}
\usepackage[margin=1.15in]{geometry}
\usepackage[colorlinks=true,linkcolor=blue,citecolor=blue,urlcolor=blue]{hyperref}
\hypersetup{
  pdftitle={On a conjecture of Fiedorowicz},
  pdfauthor={Oren Ben-Bassat},
  pdfsubject={Classifying spaces of finite aperiodic monoids},
  pdfkeywords={classifying spaces, finite monoids, aperiodic monoids,
    finite CW complexes, mapping cones, Fiedorowicz conjecture}
}

\newtheorem{theorem}{Theorem}[section]
\newtheorem{lemma}[theorem]{Lemma}
\newtheorem{proposition}[theorem]{Proposition}
\newtheorem{corollary}[theorem]{Corollary}
\newtheorem{conjecture}[theorem]{Conjecture}

\newtheorem*{theoremA}{Theorem A}
\newtheorem*{theoremB}{Theorem B}
\newtheorem*{theoremC}{Theorem C}
\theoremstyle{definition}
\newtheorem{definition}[theorem]{Definition}
\newtheorem{example}[theorem]{Example}
\theoremstyle{remark}
\newtheorem{remark}[theorem]{Remark}

\newcommand{\CP}{\mathbb{CP}}
\newcommand{\RP}{\mathbb{RP}}
\newcommand{\A}{\mathcal{A}}
\newcommand{\B}{\mathcal{B}}
\newcommand{\C}{\mathcal{C}}
\newcommand{\El}{\mathrm{El}}
\newcommand{\Cof}{\mathrm{Cof}}
\newcommand{\hol}{\mathrm{hol}}
\newcommand{\op}{\mathrm{op}}
\newcommand{\sk}{\mathrm{sk}}
\newcommand{\Sd}{\mathrm{Sd}}
\newcommand{\id}{\mathrm{id}}
\newcommand{\Dn}{\Delta_{\le n}}
\newcommand{\pair}[2]{\langle #1,#2\rangle}
\newcommand{\ThmA}{\hyperref[thm:main-A]{Theorem~A}}
\newcommand{\ThmB}{\hyperref[thm:main-B]{Theorem~B}}
\newcommand{\ThmC}{\hyperref[thm:main-C]{Theorem~C}}

\title[On a conjecture of Fiedorowicz]{On a conjecture of Fiedorowicz}
\author{Oren Ben-Bassat}
\address{\hspace*{-\parindent}{\normalfont Oren Ben-Bassat}\newline
  Department of Mathematics\newline
  University of Haifa\newline
  Israel}

\email{\href{mailto:ben-bassat@math.haifa.ac.il}{ben-bassat@math.haifa.ac.il}}
\date{\today}

\begin{document}
\begin{abstract}
Fiedorowicz conjectured that every finite simply connected CW complex is
homotopy equivalent to the classifying space of a finite monoid.  We prove
this conjecture.  Our main construction takes a finite monoid M and a map from a
finite connected CW complex into the classifying space of M, and produces a new finite
monoid whose classifying space is homotopy equivalent to the mapping cone of
the given map.  Applying this construction successively to the attaching maps of a
finite CW model proves the result.  As an explicit example of the algorithm, we construct a
1795-element monoid whose classifying space is homotopy equivalent to two dimensional complex projective space.
\end{abstract}
\maketitle

\section{Introduction}

McDuff~\cite{McD} proved that every connected CW complex is homotopy equivalent to the
classifying space $BM$ of some discrete monoid $M$, regarded as a one-object category.
If $M$ is finite, then $\pi_1(BM)$ is finite (Section~\ref{sec:fundamental-group}).
Throughout, \emph{simply connected} means path connected with trivial fundamental group.
In 1984 Fiedorowicz~\cite[p.~323]{Fie84} proposed the following conjecture:

\begin{conjecture}\label{conj:fiedorowicz}
Every finite simply connected CW complex is homotopy equivalent to $BM$ for some finite monoid $M$.
\end{conjecture}
We prove a slightly stronger statement, \ThmC{} below.

A finite monoid is \emph{aperiodic} if every element has stationary powers: for every $a$ there is an $r\geq1$ such
that $a^r=a^{r+1}$.  Equivalently, every subgroup of $M$ is trivial; this is standard~\cite{RS}.
The monoids constructed below are aperiodic, so the result is stronger than the original
conjecture.  Fiedorowicz gave the following five-element example.

\begin{example}\label{ex:intro-sphere}
Let $L_0=\{1\}\sqcup\{x_{ij}\mid i,j\in\{0,1\}\}$, with multiplication
$x_{ij}x_{k\ell}=x_{i\ell}$ and with $1$ acting as the identity.  The four nonidentity elements form a
$2\times2$ rectangular band.  This monoid
is aperiodic and $BL_0\simeq S^2$~\cite{Fie84,Fie02}.  We recall the computation, which is used in Example~\ref{ex:filling-sphere}.  For
$P_i=\{x_{i0},x_{i1}\}$, $i\in\{0,1\}$, the right module
$\mathbb Z[P_i]=x_{ii}\mathbb Z[L_0]$ is projective because it is the image of an idempotent endomorphism.
Fiedorowicz exhibits the exact projective resolution~\cite[p.~34]{Fie02}
\[
 0\longrightarrow \mathbb Z[P_0]\oplus\mathbb Z[P_1]
 \xrightarrow{\ \iota\ }\mathbb Z[L_0]
 \xrightarrow{\ (x_{00}-x_{01})(-)\ }\mathbb Z[P_0]
 \xrightarrow{\ \epsilon\ }\mathbb Z\longrightarrow0,
\]
where $\iota$ includes the two rows and $\epsilon$ is augmentation.  Each $P_i$ is transitive, so tensoring with the
trivial left module gives
$0\to\mathbb Z^2\xrightarrow{(u,v)\mapsto u+v}\mathbb Z\xrightarrow{0}\mathbb Z\to0$, with
$\mathbb Z^2$ in degree $2$.  Hence
$H_q(BL_0;\mathbb Z)\cong\operatorname{Tor}^{\mathbb Z[L_0]}_q(\mathbb Z,\mathbb Z)$ is $\mathbb Z$ for
$q=0,2$ and is zero otherwise.
Write $U(M)$ for the universal group of a monoid $M$, namely the initial group receiving a monoid homomorphism from
$M$.  Every $x_{ij}$ is idempotent, so every homomorphism from $L_0$ to a group sends it to the identity.  Thus
$U(L_0)$, and hence $\pi_1(BL_0)$~\cite[Lemma~1]{McD}, is trivial.  The Hurewicz isomorphism now supplies a map
$S^2\to BL_0$ inducing a generator of $H_2$.  This map is an integral homology equivalence between simply connected CW
complexes, so the homological Whitehead theorem proves $BL_0\simeq S^2$.
\end{example}

Sweeney~\cite{Swe} showed that the class of spaces $BM$, with $M$ finite, is closed under suspension: his
construction $J^{\{1,2\}}(M)$ has $3|M|$ elements and satisfies $B J^{\{1,2\}}(M)\simeq\Sigma BM$.  Martinian
and Steinberg~\cite{MS} realized all finite wedges of simply connected Moore spaces by finite monoids.

Let $\Delta$ be the category of the ordered sets $[k]=\{0<1<\cdots<k\}$ and weakly monotone maps.  Write
$\Dn\subseteq\Delta$ for the full subcategory on $[0],\dots,[n]$, and let
$O_{n+1}=\Delta([n],[n])$, a monoid with $\binom{2n+1}{n}$ elements.  For a map
$\theta\colon[j]\to[k]$ in $\Delta$, we write
$\theta^{\op}\colon[k]\to[j]$ for the corresponding arrow in the opposite category $\Delta^{\op}$.

If $e$ is an idempotent of a monoid $M$, then the local monoid $eMe$ has identity $e$.  We write
$j_e\colon B(eMe)\to BM$ for its canonical morphism in the homotopy category.  Concretely, the Karoubi envelope $\operatorname{Kar}(M)$ has
the idempotents of $M$ as objects and $\operatorname{Hom}(e,e')=e'Me$.  The two morphisms represented by $e$,
one from $e$ to $1$ and one from $1$ to $e$, have composite $e=\id_e$; hence every object $e$ is a retract of $1$.
Lemma~\ref{lem:gen} shows that the inclusion of $1$ induces an
equivalence $BM\xrightarrow{\simeq}B\operatorname{Kar}(M)$, and $j_e$ is represented by the span
\[
 B(eMe)\longrightarrow B\operatorname{Kar}(M)\xleftarrow{\simeq}BM.
\]

\begin{theoremA}
\phantomsection\label{thm:main-A}%
Let $L$ be a finite monoid, let $Y$ be a finite reduced simplicial set, put $n=\dim Y$, and let
$f\colon Y\to N(L)$ be a simplicial map.  Write $*_n\in Y_n$ for the totally degenerate $n$-simplex.  There is a finite monoid
\[ M(Y,f)\;=\;O_{n+1}^{\op}\ \sqcup\ \bigl(Y_n\times L\times\{0,\dots,n\}\bigr), \]
described explicitly in Definition~\ref{def:M}.  For every $i\in[n]$, the element $e_i=(*_n,1,i)$ is idempotent, and
\[
 \lambda_i\colon L\xrightarrow{\ \cong\ }e_iM(Y,f)e_i,
 \qquad
 \lambda_i(l)=(*_n,l,i),
\]
is an isomorphism of monoids.  There is an equivalence, independent of $i$,
\[
 \Phi\colon BM(Y,f)\xrightarrow{\ \simeq\ }
 \Cof\bigl(|f|\colon |Y|\to BL\bigr)
\]
such that, for every $i\in[n]$,
\[
 \Phi\circ j_{e_i}\circ B\lambda_i\simeq\iota_f\colon BL\longrightarrow
 \Cof\bigl(|f|\colon |Y|\to BL\bigr),
\]
where $\iota_f$ is the canonical map from $BL$ to the homotopy cofiber.  Thus $\Phi$ is an equivalence under $BL$.
\end{theoremA}

In words, the left action of $O_{n+1}^{\op}$ on $Y_n\times L\times\{0,\dots,n\}$ pulls a top simplex back along a
simplicial operator and transports its label $l$ by the $f$-holonomy of the corresponding edge path
(Section~\ref{sec:monoid}); Examples~\ref{ex:holonomy} and~\ref{ex:filling-sphere} show this in the smallest cases.

Note that $\lambda_i$ sends $1$ to $e_i$ rather than to the identity of $M(Y,f)$; this is why the map
$BL\to BM(Y,f)$ passes through $B(e_iM(Y,f)e_i)$.  Since $Y_n$ includes the degenerate simplices,
\[
 \#M(Y,f)=\binom{2n+1}{n}+(n+1)(\#L)(\#Y_n).
\]

\begin{theoremB}
\phantomsection\label{thm:main-B}%
Let $L$ be a finite aperiodic monoid, let $Y$ be a finite connected CW complex, and let
$f\colon Y\to BL$ be a map.  One can choose a finite reduced simplicial set $S(Y,f)$, a simplicial map
$f_S\colon S(Y,f)\to N(L)$, and a homotopy equivalence
$h_{Y,f}\colon|S(Y,f)|\to Y$ such that $|f_S|\simeq f\circ h_{Y,f}$.  The monoid
$M(S(Y,f),f_S)$ is finite and aperiodic, and
\[
 BM(S(Y,f),f_S)\simeq\Cof(|f_S|)\simeq\Cof(f).
\]
The symbols $S(Y,f)$, $f_S$, and $h_{Y,f}$ denote choices; no functoriality in $(Y,f)$ is asserted.
\end{theoremB}

\begin{theoremC}
\phantomsection\label{thm:main-C}%
Conjecture~\ref{conj:fiedorowicz} holds.  More strongly, every finite simply connected CW complex is homotopy
equivalent to $BM$ for a finite aperiodic monoid $M$.
\end{theoremC}

\ThmC{} follows from \ThmB{} by induction over cells, starting from the trivial monoid $\{1\}$.  Taking $L=\{1\}$ in \ThmA{}
gives:

\begin{corollary}\label{cor:susp}
If $Y$ is a finite reduced simplicial set of dimension $n$, let
\[
 M_Y=O_{n+1}^{\op}\sqcup\bigl(Y_n\times\{0,\dots,n\}\bigr).
\]
The multiplication on $M_Y$ is determined by
\[
\begin{aligned}
 \theta^{\op}\phi^{\op}&=(\phi\circ\theta)^{\op}, &
 \theta^{\op}(y,i)&=(\theta^*y,i),\\
 (y,i)\theta^{\op}&=(y,\theta(i)), &
 (y,i)(y',i')&=(y,i'),
\end{aligned}
\]
for $\theta,\phi\in O_{n+1}$, $y,y'\in Y_n$, and $i,i'\in[n]$.  The reversed order in the first formula is
composition in the opposite category, and the identity is $\id_{[n]}^{\op}$.  Then
$BM_Y\simeq\Sigma|Y|$.  Consequently, the suspension of every finite connected CW complex is homotopy equivalent
to $BM$ for a finite aperiodic monoid $M$.
\end{corollary}
\begin{proof}
The first assertion is \ThmA{} with $L=\{1\}$.  For the second, let $X$ be a finite connected CW complex.  Choose a
finite connected simplicial set $K$ and a homotopy
equivalence $|K|\simeq X$~\cite[Theorem~2C.5]{Hatcher}.  Apply Lemma~\ref{lem:reduced} with $L=\{1\}$ to the unique simplicial map
$K\to N(\{1\})$.  The lemma produces a finite reduced simplicial set $Y$ with $|Y|\simeq|K|$, and hence
$BM_Y\simeq\Sigma|Y|\simeq\Sigma X$.  The monoid $M_Y$ is aperiodic by Lemma~\ref{lem:aperiodic}.
\end{proof}

This corollary implies the theorem of Martinian and Steinberg~\cite{MS}, since a finite wedge of simply connected Moore spaces of finitely
generated abelian groups is the suspension of a finite connected complex.

\begin{remark}\label{rem:derived-group-completion}
Burke proves $|\mathbb L U(M)|\simeq\Omega BM$ for the derived universal-group completion of a discrete
monoid~\cite[Proposition~3.84]{Burke}.  Since stationary powers force $U(M)=1$ for an aperiodic monoid, \ThmC{}
implies that every based finite simply connected CW complex $X$ admits a finite monoid $M$ with
$U(M)=1$ but $|\mathbb L U(M)|\simeq\Omega X$.  Derived, rather than ordinary, group completion retains the
topology constructed here.
\end{remark}

\medskip\noindent\textbf{Outline.} Section~\ref{sec:prelim} collects facts about classifying
spaces of categories and homotopy cofibers. Section~\ref{sec:monoid} constructs $M(Y,f)$ from
the holonomy of $f$, and Section~\ref{sec:proofA} proves \ThmA. Section~\ref{sec:cells} proves
\ThmB{} and \ThmC. Section~\ref{sec:cellular-examples} gives examples, including a
$1795$-element monoid with $BM\simeq\CP^2$, and Section~\ref{sec:fundamental-group} discusses
the case of nontrivial finite fundamental group.

\section{Preliminaries}\label{sec:prelim}

\medskip
\noindent\textbf{Conventions.} A simplicial set is \emph{finite} if it has finitely many nondegenerate simplices, and it
is \emph{reduced} if it has exactly one vertex.  For a small category
$\C$, $B\C=|N\C|$, and $f\circ g$ means ``first $g$, then $f$.''  We write $\simeq$ for homotopy equivalence;
between spaces, $\cong$ denotes a homeomorphism, and elsewhere it denotes an isomorphism.

For a monoid $L$, let $\mathbf L$ be the one-object category with unique object $\bullet_L$ and with
$l\circ l'=ll'$.  We write $N(L)$ for the bar simplicial set with $N(L)_k=L^k$, with $d_0$ and $d_k$
deleting the first and last entry,
$d_i(l_1,\dots,l_k)=(\dots,l_i l_{i+1},\dots)$ for $0<i<k$, and degeneracies inserting $1$.  Thus
$N(L)=N(\mathbf L^{\op})$.  Reversal of arrows and of the affine coordinates of each simplex gives a canonical
homeomorphism $\operatorname{rev}_L:B\mathbf L\xrightarrow{\cong}|N(L)|$.  Whenever we write $BL$, or compare a
functor to $\mathbf L$ with a map to $N(L)$, we use $\operatorname{rev}_L$.  This convention removes any implicit
reversal of multiplication.  A simplicial map
$f\colon Y\to N(L)$ from a reduced $Y$ amounts to a labelling of the edges of $Y$ by elements of $L$ such that the
degenerate edge is labelled by $1$ and
$f(d_1\sigma)=f(d_2\sigma)f(d_0\sigma)$ for every $2$-simplex $\sigma$.

We use Quillen's Theorem~A~\cite{Qui} in the
form: a functor $F\colon\C\to\mathcal D$ induces a homotopy equivalence $B\C\to B\mathcal D$ if the comma category
$(F\downarrow d)$ is contractible for every object $d$, or if $(d\downarrow F)$ is contractible for every $d$.  A natural
transformation between functors induces a homotopy between the maps of classifying spaces; in particular a category with
an initial or terminal object is contractible.  For a map $f\colon X\to Z$, we write
$\Cof(f):=Z\mathop{\cup}\limits_X^h *$ for its homotopy cofiber.  When $X$ and $Z$ are CW complexes, this homotopy pushout is represented by the ordinary
mapping cone $Z\cup_f CX$.  Indeed, the cone inclusion $X\hookrightarrow CX$ is a cofibration and $CX\to *$ is a homotopy equivalence.  Hence the ordinary pushout of $Z\leftarrow X\hookrightarrow CX$ computes
the homotopy pushout of $Z\leftarrow X\to *$.  If $X$ and $Z$ merely have CW type, choose CW complexes
$X_{\mathrm{cw}}$ and $Z_{\mathrm{cw}}$ and homotopy equivalences $u\colon X_{\mathrm{cw}}\to X$ and
$v\colon Z_{\mathrm{cw}}\to Z$.  Given a homotopy inverse $w\colon Z\to Z_{\mathrm{cw}}$ of $v$, the map
$f_{\mathrm{cw}}:=w\circ f\circ u$ satisfies $v\circ f_{\mathrm{cw}}\simeq f\circ u$.
Lemma~\ref{lem:cofiber-invariance} gives
$\Cof(f_{\mathrm{cw}})\simeq\Cof(f)$.

\begin{lemma}\label{lem:cofiber-invariance}
Suppose that the square
\[
\begin{array}{ccc}
 X&\xrightarrow{f}&Z\\
 {\scriptstyle u}\big\downarrow&&\big\downarrow{\scriptstyle v}\\
 X'&\xrightarrow{f'}&Z'
\end{array}
\]
commutes up to a specified homotopy and that $u$ and $v$ are homotopy equivalences.  The square and the specified
homotopy induce a homotopy equivalence $\Cof(f)\simeq\Cof(f')$.
\end{lemma}
\begin{proof}
This is homotopy invariance of the double mapping cylinder, a model for the homotopy pushout; see
\cite[pp.~80--81]{May}.  The specified homotopy supplies the map between the two double mapping cylinders, and its
three object maps $u$, $v$, and $\id_*$ are homotopy equivalences.
\end{proof}

\begin{lemma}\label{lem:gluing}
Let $i\colon A\to X$ and $g\colon A\to B$ be maps between CW complexes, and let $P=X\cup_A B$ be their pushout.
Suppose that $g$ is a homotopy equivalence and that either $i$ or $g$ is a cofibration.  Then the canonical map
$X\to P$ is a homotopy equivalence.
\end{lemma}
\begin{proof}
The double mapping cylinder represents the homotopy pushout, and its quotient to the ordinary pushout is a homotopy
equivalence when either leg is a cofibration~\cite[pp.~80--81]{May}.  Since $g$ is a homotopy equivalence, homotopy
invariance identifies that homotopy pushout with $X$.
\end{proof}

\begin{definition}\label{def:collage}
Let $\A,\B$ be small categories and let $\Pi\colon\A^{\op}\times\B\to\mathrm{Set}$ be a functor.  The
\emph{collage} $\C(\Pi)$ of $\Pi$ has object set $\mathrm{Ob}\,\A\sqcup\mathrm{Ob}\,\B$.  For
$a,a'\in\mathrm{Ob}\,\A$ and $b,b'\in\mathrm{Ob}\,\B$, all of its morphism sets are
\[
 \C(\Pi)(a,a')=\A(a,a'),\quad
 \C(\Pi)(b,b')=\B(b,b'),\quad
 \C(\Pi)(a,b)=\Pi(a,b),\quad
 \C(\Pi)(b,a)=\varnothing.
\]
Composition within $\A$ and within $\B$ is unchanged.  If $x\in\Pi(a,b)$,
$\alpha\colon a'\to a$ is a morphism of $\A$, and $\beta\colon b\to b'$ is a morphism of $\B$, then
\[
 x\circ\alpha:=\Pi(\alpha,\id_b)(x),\qquad
 \beta\circ x:=\Pi(\id_a,\beta)(x),\qquad
 \beta\circ x\circ\alpha=\Pi(\alpha,\beta)(x).
\]
These formulas, together with the two internal composition laws, define every possible composition in $\C(\Pi)$;
the category axioms follow from functoriality of $\Pi$.

The \emph{category of elements} $\El(\Pi)$ has objects $(a,x,b)$ with
$a\in\mathrm{Ob}\,\A$, $b\in\mathrm{Ob}\,\B$, and $x\in\Pi(a,b)$.  A morphism
$(a,x,b)\to(a',x',b')$ is a pair
$(\alpha,\beta)\in\A(a,a')\times\B(b,b')$ satisfying
$\Pi(\alpha,\id_{b'})(x')=\Pi(\id_a,\beta)(x)$.  Equivalently, $x'\circ\alpha=\beta\circ x$ in $\C(\Pi)$.
The projection functors $s\colon\El(\Pi)\to\A$ and $t\colon\El(\Pi)\to\B$ send $(a,x,b)$ to $a$ and $b$,
respectively, and send $(\alpha,\beta)$ to $\alpha$ and $\beta$.
\end{definition}

\begin{lemma}\label{lem:gen}
Let $\mathcal D$ be a small category and $c_0$ an object such that every object of $\mathcal D$ is a retract of
$c_0$.  Then the inclusion of the monoid $\mathcal D(c_0,c_0)$, viewed as a category with one object, induces a
homotopy equivalence
$B\,\mathcal D(c_0,c_0)\simeq B\mathcal D$.
\end{lemma}
\begin{proof}
Let $\iota$ denote the inclusion of the one-object category $\mathcal D(c_0,c_0)$ into $\mathcal D$.  To apply
Quillen's Theorem~A, fix an object $c$ of $\mathcal D$ and consider the comma category $(\iota\downarrow c)$.  The
objects of this comma category are the
arrows $\varphi\colon c_0\to c$; a morphism from $\varphi$ to $\varphi'$ is an endomorphism
$m\colon c_0\to c_0$ such that $\varphi'\circ m=\varphi$.

When $c=c_0$, the object $\id_{c_0}$ is terminal, so $(\iota\downarrow c_0)$ is contractible.  For a general $c$,
choose arrows $s\colon c\to c_0$ and $r\colon c_0\to c$ with $r\circ s=\id_c$.  Postcomposition gives functors
\[
 (\iota\downarrow c)\xrightarrow{\ s\circ(\ )\ }(\iota\downarrow c_0)
 \xrightarrow{\ r\circ(\ )\ }(\iota\downarrow c)
\]
whose composite is the identity.  Hence $B(\iota\downarrow c)$ is a retract of the contractible space
$B(\iota\downarrow c_0)$.  The identity map of $B(\iota\downarrow c)$ factors through a contractible space
and is nullhomotopic, so
$B(\iota\downarrow c)$ is itself contractible.  By Quillen's Theorem~A~\cite{Qui}, $B\iota$ is a homotopy
equivalence.
\end{proof}

\begin{lemma}\label{lem:collage}
Let $\A$ and $\B$ be small categories, let $\Pi\colon\A^{\op}\times\B\to\mathrm{Set}$ be a functor, and let
$\C(\Pi)$ and $\El(\Pi)$ be its collage and category of elements from Definition~\ref{def:collage}.  If
$s\colon\El(\Pi)\to\A$ and $t\colon\El(\Pi)\to\B$ are the projection functors, then $B\C(\Pi)$ is homotopy
equivalent to the homotopy pushout
\[
 B\A\xleftarrow{\ Bs\ }B\El(\Pi)\xrightarrow{\ Bt\ }B\B,
\]
compatibly with the inclusions $B\A\to B\C(\Pi)$ and $B\B\to B\C(\Pi)$.
\end{lemma}
\begin{proof}
Let $\operatorname{Cyl}(s,t)$ be the Grothendieck construction of the span
$\A\xleftarrow{s}\El(\Pi)\xrightarrow{t}\B$.  Thomason's homotopy-colimit theorem gives, naturally in the span,
\[
 B\operatorname{Cyl}(s,t)\simeq
 \operatorname*{hocolim}\bigl(B\A\xleftarrow{Bs}B\El(\Pi)\xrightarrow{Bt}B\B\bigr)
 \qquad\text{\cite[Theorem~1.2]{Tho}}.
\]
In particular, this equivalence respects the maps from $B\A$ and $B\B$.

There is a comparison functor $F\colon\operatorname{Cyl}(s,t)\to\C(\Pi)$.  The functor $F$ is the identity on $\A$
and $\B$, and sends an object
$(a,x,b)$ of $\El(\Pi)$ to $a$.  On an arrow $(\alpha,\beta)$ in $\El(\Pi)$ it takes the $\A$-component $\alpha$;
an arrow $\alpha\in\operatorname{Cyl}(s,t)(e,a')$ is sent to the same arrow in $\A$, while
$\beta\in\operatorname{Cyl}(s,t)((a,x,b),b')$ is sent to $\beta\circ x\colon a\to b'$.  The defining equation
$x'\circ\alpha=\beta\circ x$ in $\El(\Pi)$ proves that $F$ respects the mixed compositions.

If $c=b_0\in\B$, then $(b_0\downarrow F)=(b_0\downarrow\B)$ has the initial object $(b_0,\id_{b_0})$.
Now let $c=a_0\in\A$.  The objects of
$(a_0\downarrow F)$ are of three kinds:
$(a,\alpha\colon a_0\to a)$, $((a,x,b),\alpha\colon a_0\to a)$, and $(b,y)$ with $y\in\Pi(a_0,b)$.  Denote the corresponding full subcategories by
$\A_{a_0},\mathcal E_{a_0},\B_{a_0}$, respectively.  Then
$(a_0\downarrow F)$ is the Grothendieck construction of the span
$\A_{a_0}\leftarrow\mathcal E_{a_0}\to\B_{a_0}$: all mixed morphisms start in $\mathcal E_{a_0}$, and its two
functors are $((a,x,b),\alpha)\mapsto(a,\alpha)$ and
$p_{a_0}\colon((a,x,b),\alpha)\mapsto(b,x\circ\alpha)$.  Now
$\A_{a_0}=(a_0\downarrow\A)$ has the initial object $(a_0,\id_{a_0})$, and $p_{a_0}$ is a homotopy equivalence: the section
$\sigma_{a_0}(b,y)=((a_0,y,b),\id)$ satisfies $p_{a_0}\circ\sigma_{a_0}=\id$, and the arrows
$(\alpha,\id_b)$ give $\sigma_{a_0}\circ p_{a_0}\Rightarrow\id$.  Thomason's theorem identifies
$B(a_0\downarrow F)$ with the homotopy pushout of
$B\A_{a_0}\leftarrow B\mathcal E_{a_0}\to B\B_{a_0}$.  Since $B\A_{a_0}$ is contractible and
$Bp_{a_0}$ is a homotopy equivalence, this homotopy pushout is contractible.  Thus every comma category required by
Quillen's Theorem~A is contractible.  It follows that $BF$ is a homotopy equivalence~\cite{Qui}, and the natural
homotopy-colimit equivalence above proves the claim, including compatibility with the two structure maps.
\end{proof}

\begin{lemma}\label{lem:trunc}
Let $Y$ be a simplicial set all of whose nondegenerate simplices have dimension $\le n$, and let $\Delta\!\downarrow\! Y$
be its category of simplices, whose objects are $([k],y)$, with $y\in Y_k$, and whose morphisms are the maps
\[
 \theta\colon([j],\theta^*y)\longrightarrow([k],y).
\]
Then the
inclusion $\Dn\!\downarrow\! Y\to\Delta\!\downarrow\! Y$ is a homotopy equivalence.  Consequently the restriction
of the last-vertex map gives
\[
 \operatorname{lv}_n:=\left.\operatorname{lv}\right|_{N(\Dn\downarrow Y)}
 \colon N(\Dn\!\downarrow\!Y)\longrightarrow Y,
 \qquad B(\Dn\!\downarrow\! Y)\simeq|Y|,
\]
where $\operatorname{lv}\colon N(\Delta\!\downarrow\! Y)\to Y$ is the full last-vertex map.
\end{lemma}
\begin{proof}
By Quillen's Theorem~A it suffices to show that, for each simplex $([m],\sigma)$ of $Y$, the following category of
its factorizations through simplices of dimension at most $n$ is contractible.  Let $\mathcal F$ be the category of triples
$(k,\tau,\theta)$ with $k\le n$, $\tau\in Y_k$,
$\theta\colon[m]\to[k]$ and $\theta^*\tau=\sigma$; a morphism $(k,\tau,\theta)\to(k',\tau',\theta')$ is an
$\alpha\colon[k]\to[k']$ with $\alpha^*\tau'=\tau$ and $\alpha\circ\theta=\theta'$.  This is the comma category
$(([m],\sigma)\downarrow\iota)$ for the inclusion
$\iota\colon\Dn\!\downarrow\!Y\to\Delta\!\downarrow\!Y$.  Quillen's
Theorem~A~\cite{Qui} reduces the desired homotopy equivalence to proving that $\mathcal F$ is contractible.

Let $\mathcal F_s\subseteq\mathcal F$ be the full subcategory spanned by the objects
$(k,\tau,\theta)$ for which $\theta\colon[m]\to[k]$ is surjective.  Factor $\theta$ uniquely as
\[
 [m]\xrightarrow{\ \rho\ }[\bar k]\xrightarrow{\ \delta\ }[k],
\]
where $\bar k=|\operatorname{im}(\theta)|-1$, $\rho$ is surjective, and $\delta$ is the injective map enumerating
$\operatorname{im}(\theta)$ in its induced order.  The assignment
$(k,\tau,\theta)\mapsto(\bar k,\delta^*\tau,\rho)$ extends to a retraction
$r\colon\mathcal F\to\mathcal F_s$.  A morphism $\alpha$ carries
$\operatorname{im}(\theta)$ into $\operatorname{im}(\theta')$ and induces a unique weakly monotone map
$\bar\alpha$ between the two ordered images.  The identities
$\bar\alpha\circ\rho=\rho'$, $\delta'\circ\bar\alpha=\alpha\circ\delta$, and
$\bar\alpha^*\delta'^*\tau'=\delta^*\tau$ show that $r$ is functorial.  The injections
$\delta$ define a natural transformation
$\iota\circ r\Rightarrow\id_{\mathcal F}$, where $\iota\colon\mathcal F_s\hookrightarrow\mathcal F$ is the inclusion.
Thus $B\mathcal F\simeq B\mathcal F_s$.
The Eilenberg--Zilber lemma~\cite[pp.~26--27]{GZ} says that $\sigma$ is a unique degeneracy of a unique
nondegenerate simplex: $\sigma=\pi^*\sigma_0$ for a unique surjection $\pi\colon[m]\to[k_0]$ and a unique
nondegenerate $\sigma_0$.  Since $Y$ has dimension at most $n$, we have $k_0\le n$.  If
$(k,\tau,\rho)\in\mathcal F_s$ and $\tau=\pi'^*\tau_0$ with $\tau_0$ nondegenerate, then
$\sigma=(\pi'\circ\rho)^*\tau_0$.  Uniqueness gives $\tau_0=\sigma_0$ and $\pi'\circ\rho=\pi$; surjectivity of
$\rho$ makes $\pi'$ the unique morphism to $(k_0,\sigma_0,\pi)$.  Hence this object is terminal in $\mathcal F_s$.

Let $\iota\colon\Dn\!\downarrow\!Y\to\Delta\!\downarrow\!Y$ denote the inclusion just shown to induce a homotopy
equivalence on classifying spaces.  The full last-vertex map
$\operatorname{lv}\colon N(\Delta\!\downarrow\!Y)\to Y$ is a weak equivalence by
\cite[Proposition~4.2]{Tho}, and the restricted last-vertex map satisfies the strict identity
$\operatorname{lv}_n=\operatorname{lv}\circ N\iota$.
Since $N\iota$ and $\operatorname{lv}$ are weak equivalences, their composite $\operatorname{lv}_n$ is a weak
equivalence.  After geometric realization, its source and target are CW complexes, so Whitehead's theorem gives
$B(\Dn\!\downarrow\!Y)\simeq|Y|$.
\end{proof}

\section{The monoid \texorpdfstring{$M(Y,f)$}{M(Y,f)}}\label{sec:monoid}

Fix $L$, $Y$, $f$ and $n$ as in \ThmA.  For $y\in Y_k$ write $f(y)=(f_1(y),\dots,f_k(y))\in L^k$ and put
$\hol_i(y)=f_{i+1}(y)\cdots f_k(y)$ for $0\le i\le k$; this is the product of the labels of the edges of $y$ from
vertex $i$ to the last vertex.  If
$\theta\colon[j]\to[k]$, then, for $1\le r\le j$, the $r$th edge of $\theta^*y$ has label
$f_{\theta(r-1)+1}(y)\cdots f_{\theta(r)}(y)$, with an empty product interpreted as $1$.  Multiplying these edge labels gives, for
$\psi\colon[i]\to[j]$,
\begin{equation}\label{eq:cocycle}
\hol_{\psi(i)}(\theta^*y)\cdot\hol_{\theta(j)}(y)=\hol_{\theta(\psi(i))}(y).
\end{equation}
Equation~\eqref{eq:cocycle} says that the edge-path from the vertex $\theta(\psi(i))$ to the last vertex of $y$ may be split
at the vertex $\theta(j)$.

Define a simplicial right $L$-set $E_f$ degreewise by $(E_f)_k=Y_k\times L$.  Its right action is
$(y,l)\cdot m=(y,lm)$, and for $\theta\colon[j]\to[k]$ its simplicial operator is
$\theta^*(y,l)=\bigl(\theta^*y,\hol_{\theta(j)}(y)l\bigr)$.
Identity~\eqref{eq:cocycle} verifies the simplicial identities.  Equivalently,
$E_f=Y\times_{N(L)}EL$, the pullback of the standard universal bar right $L$-set $EL\to N(L)$.  Each
$(E_f)_k=Y_k\times L$ is a disjoint union of copies of $L$ acting on itself by right multiplication, and
$(E_f)_k/L=Y_k$.  Since $Y$ is reduced, $(E_f)_0=\{*\}\times L$.  For $(y,l)\in(E_f)_k$ and $0\le i\le k$,
define $\pair{i}{(y,l)}\in L$ by $i^*(y,l)=(*,\pair{i}{(y,l)})$; equivalently,
$\pair{i}{(y,l)}=\hol_i(y)l$,
where $i$ also denotes the map $[0]\to[k]$ with value $i$.  Then
$\pair{\theta(i)}{(y,l)}=\pair{i}{\theta^*(y,l)}$ and
$\pair{i}{(y,l)\cdot m}=\pair{i}{(y,l)}m$.

\begin{definition}\label{def:retraction-category}
The right $L$-action and simplicial operators make $E_f$ a functor
\[
 E_f\colon\mathbf L^{\op}\times\Dn^{\op}\longrightarrow\mathrm{Set}.
\]
Write $\mathcal C(Y,f):=\C(E_f)$ for the collage of $E_f$.
The objects of $\mathcal C(Y,f)$ are $\bullet_L,[0],\ldots,[n]$; the full subcategories on
$\bullet_L$ and on $[0],\ldots,[n]$ are respectively $\mathbf L$ and $\Dn^{\op}$; and the only remaining arrows are
$(y,l)\colon\bullet_L\to[k]$ for $(y,l)\in(E_f)_k$.  The collage composition is
$(y,l)\circ m=(y,l)\cdot m$ and $\theta^{\op}\circ(y,l)=\theta^*(y,l)$.
Let $s_0=(*,1)\colon\bullet_L\to[0]$.  Define $\mathcal C^{\mathrm{ret}}(Y,f)$ by adjoining to
$\mathcal C(Y,f)$ one arrow $\rho\colon[0]\to\bullet_L$
subject only to the relation $\rho\circ s_0=\id_{\bullet_L}$.  The arrow $\rho$ makes $\bullet_L$ a retract of $[0]$.
\end{definition}

Lemma~\ref{lem:collage} identifies $B\mathcal C(Y,f)$ with the homotopy pushout
$BL\leftarrow B\El(E_f)\to B\Dn^{\op}$.  To pass from the category $\mathcal C(Y,f)$ to a monoid, we adjoin a
retraction.  Lemma~\ref{lem:frame-splitting} will show that adjoining the retraction $\rho$ does not change its
classifying-space homotopy type.  Lemma~\ref{lem:retract} will then replace
$\mathcal C^{\mathrm{ret}}(Y,f)$ by the endomorphism monoid of $[n]$.

\begin{lemma}\label{lem:retraction-normal-forms}
Every arrow of $\mathcal C^{\mathrm{ret}}(Y,f)$ has a unique reduced representative.  The reduced representatives
that contain $\rho$ are exactly
\[
 a\circ\rho\circ i^{\op}\colon[k]\longrightarrow\bullet_L
 \quad\text{or}\quad
 (y,l)\circ\rho\circ i^{\op}\colon[k]\longrightarrow[j],
\]
where $a\in L$, $(y,l)\in(E_f)_j$, and $i\colon[0]\to[k]$ selects a vertex.  Consequently
\begin{gather*}
\mathcal C^{\mathrm{ret}}(Y,f)(\bullet_L,\bullet_L)=L,
\quad \mathcal C^{\mathrm{ret}}(Y,f)(\bullet_L,[k])=(E_f)_k,
\quad \mathcal C^{\mathrm{ret}}(Y,f)([k],\bullet_L)=\{0,\dots,k\}\times L,\\
\mathcal C^{\mathrm{ret}}(Y,f)([k],[j])
=\Delta([j],[k])\ \sqcup\ \bigl((E_f)_j\times\{0,\dots,k\}\bigr).
\end{gather*}
Under these identifications write $(i,a)=a\circ\rho\circ i^{\op}$ and
$(y,l)\star i=(y,l)\circ\rho\circ i^{\op}$.
Composition is obtained by concatenating these representatives and reducing with
$\rho\circ s_0=\id_{\bullet_L}$.
\end{lemma}
\begin{proof}
Work modulo the composition laws in the collage $\mathcal C(Y,f)$: compose every maximal string of consecutive
collage arrows, so that a word alternates between collage arrows and occurrences of $\rho$.  Every element of
$(E_f)_0$ is uniquely $s_0\cdot a$, and we orient the new relation as
$\rho\circ(s_0\cdot a)\longrightarrow a$.
The rule strictly decreases the number of occurrences of $\rho$, so every reduction sequence terminates.  The only
overlap between this rule and collage composition is resolved because $(\rho\circ s_0)\circ a=a$ and
$\rho\circ(s_0\circ a)=\rho\circ(s_0\cdot a)\longrightarrow a$.
Reductions at distinct occurrences of $\rho$ commute.  Induction on their number shows that every word
reduces to a unique irreducible word.

Two occurrences of $\rho$ would be separated by an arrow $\bullet_L\to[0]$ and would expose a reducible
subword.  Thus an irreducible word contains at most one $\rho$.  On the source side of that $\rho$, the arrow must be
a unique vertex map $i^{\op}\colon[k]\to[0]$; on the target side, the composed collage arrow is uniquely
either $a\colon\bullet_L\to\bullet_L$ or $(y,l)\colon\bullet_L\to[j]$.  The reduced words containing $\rho$ are
exactly $a\circ\rho\circ i^{\op}\colon[k]\to\bullet_L$ and
$(y,l)\circ\rho\circ i^{\op}\colon[k]\to[j]$,
where $a\in L$, $(y,l)\in(E_f)_j$, and $i\colon[0]\to[k]$.  Words without $\rho$ are the original collage arrows;
separating these arrows and the two displayed families by source and target gives the four Hom-set formulas in the
statement.  Associativity is automatic because
$\mathcal C^{\mathrm{ret}}(Y,f)$ is a category presented by these generators and relations.
\end{proof}

\begin{definition}\label{def:M}
With $n=\dim Y$, define
\[
 M(Y,f)=\operatorname{End}_{\mathcal C^{\mathrm{ret}}(Y,f)}([n]).
\]
By Lemma~\ref{lem:retraction-normal-forms}, its underlying set is
\[
 M(Y,f)=O_{n+1}^{\op}\sqcup\bigl((E_f)_n\times\{0,\dots,n\}\bigr)
 =O_{n+1}^{\op}\sqcup\bigl(Y_n\times L\times\{0,\dots,n\}\bigr).
\]
For $(y,l)\in(E_f)_n$ and $i\in[n]$, write
$(y,l,i):=(y,l)\circ\rho\circ i^{\op}\in M(Y,f)$.  The subset
$(E_f)_n\times\{0,\dots,n\}$ is a two-sided ideal.
Concatenating the normal forms of Lemma~\ref{lem:retraction-normal-forms} and using
$\rho\circ(s_0\cdot m)=m$ gives the three multiplication formulas.  For example,
$i^{\op}\circ(y',l')=i^*(y',l')=s_0\cdot(\hol_i(y')l')$, which gives the ideal-by-ideal product in the next
display.  Explicitly,
\begin{gather*}
(y,l,i)\,(y',l',i')=\bigl(y,\ l\,\hol_i(y')\,l',\ i'\bigr),\\
\theta^{\op}\,(y,l,i)=\bigl(\theta^*y,\ \hol_{\theta(n)}(y)\,l,\ i\bigr),\qquad
(y,l,i)\,\theta^{\op}=\bigl(y,\ l,\ \theta(i)\bigr).
\end{gather*}
Here $\theta\in O_{n+1}$, $y,y'\in Y_n$, $l,l'\in L$, and $i,i'\in[n]$.
The identity is $\id_{[n]}^{\op}\in O_{n+1}^{\op}$, and associativity follows because these are compositions in
$\mathcal C^{\mathrm{ret}}(Y,f)([n],[n])$.
\end{definition}

\begin{example}\label{ex:holonomy}
When $L=\{1\}$, as in Corollary~\ref{cor:susp}, every holonomy $\hol_i(y)$ is trivial.  The
smallest case in which the holonomy matters is the following.  \ThmA{} does not require $L$
to be aperiodic, so take $L=C_2=\{1,g\}$.  Let $Y=\Delta^1/\partial\Delta^1$, with nondegenerate
edge $s$, and for $a\in C_2$ let $f_a\colon Y\to N(C_2)$ send $s$ to $a$.  Then $n=1$,
$O_2=\{\id,c_0,c_1\}$ with $c_k$ constant with value $k$, and
\[
 M(Y,f_a)=\{\id^{\op},c_0^{\op},c_1^{\op}\}\sqcup\bigl(\{*_1,s\}\times C_2\times\{0,1\}\bigr)
\]
has $11$ elements for either value of $a$.  The only holonomy that can be nontrivial is
$\hol_0(s)=a$; in particular $\hol_1(s)=1$.  Since $c_0^*s=c_1^*s=*_1$,
Definition~\ref{def:M} gives
\[
 c_0^{\op}\,(s,1,i)=(*_1,a,i)=\lambda_i(a),\qquad
 c_1^{\op}\,(s,1,i)=(*_1,1,i)=e_i.
\]
Thus the two vertex operators send the lifted edge $(s,1)\in(E_{f_a})_1$ to
$\lambda_i(a)$ and $\lambda_i(1)=e_i$ in the corner $e_iM(Y,f_a)e_i\cong C_2$.

Take $a=g$.  A homomorphism $\phi$ from $M(Y,f_g)$ to a group kills the idempotents
$c_0^{\op}$, $c_1^{\op}$ and $e_i$.  The second identity then gives $\phi(s,1,i)=1$, and the
first gives $\phi(\lambda_i(g))=1$.  Since $(s,g,i)=(s,1,0)\,\lambda_i(g)$, every element maps
to $1$, so $\pi_1(BM(Y,f_g))=U(M(Y,f_g))$ is trivial~\cite[Lemma~1]{McD}.  This agrees with
\ThmA: $|f_g|$ is the inclusion $\RP^1\subseteq\RP^\infty\simeq BC_2$, and attaching a $2$-cell
along it kills $\pi_1$.  If instead $a=1$, then $c_0^{\op}(s,1,i)=e_i$, the relation between
$g$ and $1$ disappears, and $|f_1|$ is null-homotopic, so $BM(Y,f_1)\simeq BC_2\vee S^2$ has
fundamental group $C_2$.  The two monoids have the same underlying set and differ only in the
value of $\hol_0(s)$.
\end{example}

\begin{remark}\label{rem:top-simplices}
Although the underlying set in Definition~\ref{def:M} displays only $Y_n$, the monoid depends on the full simplicial
set $Y$ and on the full map $f$.  Indeed, the left action uses every pullback $\theta^*y$, while the product uses the
edge labels of all faces of $y$ through $\hol_i(y)$.  More explicitly, choose a surjection
$\sigma\colon[n]\twoheadrightarrow[k]$ and a section $\delta\colon[k]\hookrightarrow[n]$.  Then
$\sigma^*\colon Y_k\to Y_n$ is injective, and it identifies $Y_k$ with the image of the idempotent
$(\delta\circ\sigma)^*=\sigma^*\circ\delta^*$ acting on $Y_n$.  In this sense the $O_{n+1}$-action on $Y_n$ records every lower
simplicial degree and every simplicial operator; naturality of $f$ records the corresponding lower-dimensional
labels.  Taking $n=\dim Y$, as we do throughout, introduces no additional truncation choice.
\end{remark}

\begin{lemma}\label{lem:retract}
Every object of $\mathcal C^{\mathrm{ret}}(Y,f)$ is a retract of $[n]$.  Hence Lemma~\ref{lem:gen} gives
\[
 BM(Y,f)\simeq B\mathcal C^{\mathrm{ret}}(Y,f).
\]
\end{lemma}
\begin{proof}
Choose an injection $\iota\colon[k]\to[n]$ and a retraction $r\colon[n]\to[k]$ with
$r\circ\iota=\id_{[k]}$.  Then $r^{\op}$ and $\iota^{\op}$ exhibit $[k]$ as a retract of $[n]$ in
$\Delta^{\op}$.  For the remaining object $\bullet_L$, let
$s_n=(*_n,1)\colon\bullet_L\to[n]$ and
$\rho_{n,0}:=\rho\circ0^{\op}=(0,1)\colon[n]\to\bullet_L$, where $*_n$ is the totally degenerate
$n$-simplex.  These arrows satisfy
$\rho_{n,0}\circ s_n=\pair{0}{(*_n,1)}=1=\id_{\bullet_L}$.
Thus $\bullet_L$ is also a retract of $[n]$.
\end{proof}

\section{Proof of Theorem~A}\label{sec:proofA}

Lemma~\ref{lem:retract} already compares $BM(Y,f)$ with $B\mathcal C^{\mathrm{ret}}(Y,f)$.  It remains to compare
$B\mathcal C(Y,f)$ with the homotopy cofiber of $|f|$ and to show that adjoining $\rho$ induces an equivalence.

\begin{proposition}\label{prop:collage-cofiber}
There is an equivalence
\[
 \Phi_{\mathrm{col}}\colon B\mathcal C(Y,f)\xrightarrow{\ \simeq\ }
 \Cof(|f|\colon|Y|\to BL).
\]
Its composite with the inclusion $BL\to B\mathcal C(Y,f)$ is homotopic to $\iota_f$, the canonical cofiber map.
\end{proposition}
\begin{proof}
Lemma~\ref{lem:collage} identifies $B\mathcal C(Y,f)$ with the homotopy pushout
\[
 BL\longleftarrow B\El(E_f)\longrightarrow B\Dn^{\op}.
\]
Here $[0]$ is initial in $\Dn^{\op}$, so $B\Dn^{\op}$ is contractible: the unique maps from $[0]$ give a natural
contraction to that object.  The displayed homotopy pushout is the homotopy cofiber of
$B\El(E_f)\to BL$.  We must identify its source with $|Y|$ and its structure map with $|f|$.

An object of $\El(E_f)$ is a pair $([k],(y,l))$, and a morphism
$([k],(y,l))\to([j],(y',l'))$ is a pair $(m,\theta)\in L\times\Delta([j],[k])$ with
$(y',l')\cdot m=\theta^*(y,l)$.  Forgetting the $L$-coordinate defines
$\pi_Y\colon\El(E_f)\to(\Dn\!\downarrow\!Y)^{\op}$ by
$\pi_Y([k],(y,l))=([k],y)$ and $\pi_Y(m,\theta)=\theta^{\op}$.  Define a section in the other direction by
$\sigma_Y([k],y)=([k],(y,1))$ and
$\sigma_Y(\theta^{\op})=(\hol_{\theta(j)}(y),\theta)$.
Here $\theta^{\op}$ is opposite to
$\theta\colon([j],\theta^*y)\to([k],y)$ in the category of simplices.
The cocycle identity~\eqref{eq:cocycle} says that $\sigma_Y$ preserves composition.  We have
$\pi_Y\circ\sigma_Y=\id$, and the arrows
$(l,\id_{[k]})\colon([k],(y,l))\to([k],(y,1))$ give a natural transformation
$\id\Rightarrow\sigma_Y\circ\pi_Y$.  After taking classifying spaces, this
transformation is a homotopy from the identity to $B\sigma_Y\circ B\pi_Y$.  Hence $B\pi_Y$ is a homotopy
equivalence.  Reversal identifies
the classifying spaces of a category and its opposite.  Let
$\operatorname{lv}_n\colon N(\Dn\!\downarrow\!Y)\to Y$ be the last-vertex map from Lemma~\ref{lem:trunc}, and write
$\operatorname{rev}_{\Dn\downarrow Y}\colon
B((\Dn\!\downarrow\!Y)^{\op})\to B(\Dn\!\downarrow\!Y)$ for reversal.  The composite
\[
 \begin{aligned}
 E_Y\colon B\El(E_f)&\longrightarrow|Y|,\\
 E_Y&:=|\operatorname{lv}_n|\circ\operatorname{rev}_{\Dn\downarrow Y}\circ B\pi_Y
 \end{aligned}
\]
is therefore a homotopy equivalence.

We next identify the map to $BL$, including its variance.  Define
$\hol_f\colon\Dn\!\downarrow\!Y\to\mathbf L^{\op}$ by
$\hol_f(\theta\colon([j],\theta^*y)\to([k],y))=\hol_{\theta(j)}(y)$.
The functor $\hol_f$ sends every object to $\bullet_L$.  For a second morphism
represented by $\psi\colon[i]\to[j]$, functoriality in
$\mathbf L^{\op}$ amounts to identity~\eqref{eq:cocycle}.
If $\pi_L\colon\El(E_f)\to\mathbf L$ is the projection $(m,\theta)\mapsto m$, then $\sigma_Y$ satisfies the
strict identity
\[
 \pi_L\circ\sigma_Y=\hol_f^{\op},
 \qquad
 \hol_f^{\op}\colon(\Dn\!\downarrow\!Y)^{\op}\longrightarrow\mathbf L,
\]
where $\hol_f^{\op}$ is the opposite functor.  Since
$\id\Rightarrow\sigma_Y\circ\pi_Y$, it follows that
$B\pi_L\simeq B\hol_f^{\op}\circ B\pi_Y$.

On a morphism represented by
$\theta\colon[j]\to[k]$, the edge selected by $\operatorname{lv}_n$ runs in $y$ from vertex $\theta(j)$ to its last
vertex and has $f$-label $\hol_{\theta(j)}(y)$.  On a chain of simplex maps, the successive last-vertex edges have
the labels assigned by $N(\hol_f)$.  Hence
$N(\hol_f)=f\circ\operatorname{lv}_n$ in every simplicial degree.
The identities $\pi_L\circ\sigma_Y=\hol_f^{\op}$ and
$N(\hol_f)=f\circ\operatorname{lv}_n$ give a homotopy-commutative square
\[
\begin{array}{ccc}
 B\El(E_f)&\xrightarrow{\ B\pi_L\ }&B\mathbf L\vphantom{|N(L)|}\\
 {\scriptstyle E_Y\ \simeq}\big\downarrow
 &&\big\downarrow{\scriptstyle\ \cong\ \operatorname{rev}_L}\\
 |Y|&\xrightarrow{\ |f|\ }&|N(L)|.
\end{array}
\]
Consequently
$B\mathcal C(Y,f)$ is the homotopy cofiber of $|f|$.  The collage equivalence is induced by the displayed homotopy-pushout
diagram, so its restriction to $BL$ is the canonical map $\iota_f$ into that homotopy cofiber.
\end{proof}

\begin{lemma}\label{lem:frame-splitting}
Let $E\colon\Dn^{\op}\to\mathrm{Set}$ be a truncated simplicial right $L$-set, and suppose that there is an element
$s_0\in E_0$ for which the map $L\to E_0$, $l\mapsto s_0\cdot l$, is a bijection.  Regard $E$ as a functor
$\mathbf L^{\op}\times\Dn^{\op}\to\mathrm{Set}$ and let $\C(E)$ be its collage.  Form
$\mathcal C^{\mathrm{ret}}(E)$ by adjoining $\rho\colon[0]\to\bullet_L$ with
$\rho\circ s_0=\id_{\bullet_L}$.
Then the inclusion $J_E\colon\C(E)\to\mathcal C^{\mathrm{ret}}(E)$ induces a homotopy equivalence on classifying
spaces.  In particular, for $E=E_f$ the inclusion
$J_f\colon\mathcal C(Y,f)\to\mathcal C^{\mathrm{ret}}(Y,f)$ induces
$BJ_f\colon B\mathcal C(Y,f)\xrightarrow{\simeq}B\mathcal C^{\mathrm{ret}}(Y,f)$.
\end{lemma}
\begin{proof}
Use $s_0$ to identify $E_0$ with $L$, and define $\pair{i}{u}$ by
$i^*u=s_0\cdot\pair{i}{u}$.  Then
$\pair{\theta(i)}{u}=\pair{i}{\theta^*u}$ and
$\pair{i}{u\cdot m}=\pair{i}{u}m$.  Lemma~\ref{lem:retraction-normal-forms} applies verbatim and gives
$(i,l)=l\circ\rho\circ i^{\op}$ and $u\star i=u\circ\rho\circ i^{\op}$.

We shall use the following observation.  Let $\mathcal Q_L$ have object set $L$ and
\[
 \mathcal Q_L(a,a')=\{m\in L\mid a'm=a\}.
\]
It has terminal object $1$, since the unique arrow $a\to1$ is $m=a$.
If $\mathcal T\subseteq\mathcal R$ is full, $\mathcal T$ has a terminal object, and
$(\mathcal T\downarrow r)\cong\mathcal Q_L$ for every $r\notin\mathcal T$, then Quillen's Theorem~A applied to
$\mathcal T\hookrightarrow\mathcal R$ shows that $B\mathcal R$ is contractible; for $r\in\mathcal T$, the same
comma category has terminal object $(r,\id_r)$.

Quillen's Theorem~A for $J_E$ now requires two cases.  For $c=\bullet_L$, the normal forms show that the objects of
$(J_E\downarrow\bullet_L)$ are $(\bullet_L,a)$, $a\in L$, and $([k],(i,l))$.  Let
$\mathcal R_{\bullet_L}$ be the full subcategory on the former objects and
$v_l=([0],(0,l))$.  There is a retraction
\[
 \operatorname{ret}_{\bullet_L}([k],(i,l))=v_l,\qquad
 \operatorname{ret}_{\bullet_L}(u)=i^*u,\qquad
 \operatorname{ret}_{\bullet_L}(\theta^{\op})=\id,
\]
which is the identity on $\mathcal R_{\bullet_L}$.  Here the second formula concerns a mixed arrow represented by
$u\in E_k$, and the third an arrow between the $[k]$-objects.  The maps
$i^{\op}\colon([k],(i,l))\to v_l$ form a natural transformation from the identity to inclusion after retraction:
naturality amounts to $i^{\op}\circ u=i^*u$ and
$i'^{\op}\circ\theta^{\op}=(\theta(i'))^{\op}=i^{\op}$.  Thus
$B(J_E\downarrow\bullet_L)\simeq B\mathcal R_{\bullet_L}$.  In $\mathcal R_{\bullet_L}$, the full subcategory on
the $(\bullet_L,a)$ has terminal object $(\bullet_L,1)$, and its comma category over $v_l$ is $\mathcal Q_L$ via
$a\mapsto((\bullet_L,la)\to v_l)$.  The observation makes $B\mathcal R_{\bullet_L}$, and hence
$B(J_E\downarrow\bullet_L)$, contractible.

For $c=[j]$, the objects are $(\bullet_L,u)$ for $u\in E_j$,
$([k],\theta^{\op})$ for $\theta\colon[j]\to[k]$, and
$([k],u''\star i)$ for $u''\in E_j$ and $i\in[k]$.  Let $\mathcal R_j$ be full on the $(\bullet_L,u)$ together
with $t_j=([j],\id_{[j]}^{\op})$ and $v_{u''}=([0],u''\star0)$.  Retraction to $\mathcal R_j$ is given on objects by
\[
 ([k],\theta^{\op})\longmapsto t_j,\qquad
 ([k],u''\star i)\longmapsto v_{u''},
\]
and on mixed arrows by $u'\mapsto\theta^*u'$ in the first case and $u'\mapsto i^*u'$ in the second; arrows between
non-$\bullet_L$ objects go to identities.  The components $\theta^{\op}$ and $i^{\op}$ define the natural
transformation from the identity to inclusion after retraction.  Functoriality and naturality follow from
$L$-equivariance and
\[
 \theta^{\op}\circ u'=\theta^*u',\quad i^{\op}\circ u'=i^*u',\quad
 \theta'^{\op}\circ\psi^{\op}=(\psi\circ\theta')^{\op},\quad
 i'^{\op}\circ\psi^{\op}=(\psi(i'))^{\op}.
\]
Hence $B(J_E\downarrow[j])\simeq B\mathcal R_j$.  The full subcategory on the $(\bullet_L,u)$ and $t_j$ has
terminal object $t_j$, with unique arrow $u\colon(\bullet_L,u)\to t_j$, and its comma category over $v_{u''}$ is
$\mathcal Q_L$ via $a\mapsto((\bullet_L,u''\cdot a)\to v_{u''})$.  The observation again proves
$B(J_E\downarrow[j])\simeq *$.

Every $(J_E\downarrow c)$ is contractible, and Quillen's Theorem~A gives
\[
 BJ_E\colon B\C(E)\xrightarrow{\simeq}B\mathcal C^{\mathrm{ret}}(E).
\]
\end{proof}

\begin{proof}[Proof of \ThmA]
\leavevmode\phantomsection\label{proof:main-A}%
By Definition~\ref{def:M}, $M(Y,f)$ is the endomorphism monoid at $[n]$ in the retraction category.  Fix
$i\in[n]$, and let
\[
 s_n=(*_n,1)\colon\bullet_L\longrightarrow[n],\qquad
 \rho_{n,i}:=\rho\circ i^{\op}=(i,1)\colon[n]\longrightarrow\bullet_L.
\]
As in Lemma~\ref{lem:retract}, $\rho_{n,i}\circ s_n=\id_{\bullet_L}$.  Hence
$e_i=s_n\circ\rho_{n,i}=(*_n,1,i)$ is idempotent.  The standard corner isomorphism associated to this retraction is
\[
 \operatorname{End}_{\mathcal C^{\mathrm{ret}}(Y,f)}(\bullet_L)\xrightarrow{\ \cong\ }
 e_i\operatorname{End}_{\mathcal C^{\mathrm{ret}}(Y,f)}([n])e_i,
 \qquad
 l\longmapsto s_n\circ l\circ\rho_{n,i}.
\]
In the coordinates of Definition~\ref{def:M}, it is
$\lambda_i(l)=(*_n,l,i)$.  The inverse corner isomorphism sends $x$ to
$\rho_{n,i}\circ x\circ s_n$, so $\lambda_i$ identifies $L$ with the whole local monoid $e_iM(Y,f)e_i$.

Lemma~\ref{lem:retract}, Lemma~\ref{lem:frame-splitting}, and
Proposition~\ref{prop:collage-cofiber} give
\[
\begin{aligned}
 BM(Y,f)&=B\operatorname{End}_{\mathcal C^{\mathrm{ret}}(Y,f)}([n])
          \simeq B\mathcal C^{\mathrm{ret}}(Y,f),\\
 BJ_f&\colon B\mathcal C(Y,f)\xrightarrow{\simeq}B\mathcal C^{\mathrm{ret}}(Y,f),\\
 \Phi_{\mathrm{col}}&\colon B\mathcal C(Y,f)\xrightarrow{\simeq}\Cof(|f|).
\end{aligned}
\]
Let $\Phi$ be the equivalence obtained by traversing these three equivalences.  To identify the map from $BL$, use the
Karoubi envelope $\operatorname{Kar}(\mathcal C^{\mathrm{ret}}(Y,f))$.  Its morphisms from $(c,e)$ to $(d,e')$
form the corner $e'\operatorname{Hom}(c,d)e$, where the Hom-set is taken in the retraction category.
Every object of both $\mathcal C^{\mathrm{ret}}(Y,f)$ and its Karoubi envelope is a retract of $[n]$, so
Lemma~\ref{lem:gen} identifies both classifying spaces with $BM(Y,f)$.  In the Karoubi envelope,
\[
 (\bullet_L,\id_{\bullet_L})
 \underset{\rho_{n,i}}{\stackrel{s_n}{\rightleftarrows}}
 ([n],e_i)
\]
is an isomorphism because
$\rho_{n,i}\circ s_n=\id_{\bullet_L}$ and $s_n\circ\rho_{n,i}=e_i$.
The inclusion of the endomorphism monoid of $[n]$ sends the Karoubi object $e_i$ to $([n],e_i)$, while
$\lambda_i(l)=s_n\circ l\circ\rho_{n,i}$ and
$\rho_{n,i}\circ\lambda_i(l)\circ s_n=l$.
Let $\kappa_L\colon\mathbf L\hookrightarrow\mathcal C(Y,f)$ be the inclusion of the object $\bullet_L$, and regard
$J_f\circ\kappa_L$ as the corresponding inclusion into $\mathcal C^{\mathrm{ret}}(Y,f)$.  The preceding Karoubi
isomorphism gives, under $BM(Y,f)\simeq B\mathcal C^{\mathrm{ret}}(Y,f)$,
$j_{e_i}\circ B\lambda_i\simeq B(J_f\circ\kappa_L)$.
Traversing the three equivalences that define $\Phi$ therefore gives
\[
 \Phi\circ j_{e_i}\circ B\lambda_i
 \simeq\Phi_{\mathrm{col}}\circ B\kappa_L
 \simeq\iota_f,
\]
where the last homotopy is Proposition~\ref{prop:collage-cofiber}.

The equivalence $\Phi$ was defined without choosing $i$.  For every vertex $i$, the same arrow
$s_n\colon\bullet_L\to[n]$ is a section of $\rho_{n,i}$, so the displayed Karoubi-envelope isomorphism identifies
the corner determined by $e_i$ with the same object $\bullet_L$.  Thus all the maps
$j_{e_i}\circ B\lambda_i\colon BL\to BM(Y,f)$ represent the same morphism in the homotopy category, and the
displayed compatibility holds simultaneously for every $i\in[n]$.

The monoid $O_{n+1}$ is finite, and so are $L$ and $Y_n$; hence the underlying set displayed in
Definition~\ref{def:M} is finite.
\end{proof}

\section{Attaching cells}\label{sec:cells}

Lemma~\ref{lem:aperiodic} shows that the
construction preserves aperiodicity, as the cellular induction requires.  To prove \ThmB, we first represent a map from a finite connected CW complex to $BL$ by a
finite simplicial map, then replace its source by a finite reduced simplicial set without changing the realized map.
\ThmA{} then constructs the required monoid.  \ThmC{}
follows by starting with the trivial monoid and applying \ThmB{} successively to the attaching maps of the cells.

\begin{lemma}\label{lem:finite-simplicial-approximation}
Let $K$ be a finite simplicial set, let $X$ be any simplicial set, and let
$g\colon|K|\to|X|$ be continuous.  There are a finite simplicial set $K'$, a simplicial map
$\phi\colon K'\to X$, and a homotopy equivalence $v\colon|K'|\to|K|$ such that
$|\phi|\simeq g\circ v$.
\end{lemma}

\begin{proof}
Apply relative simplicial approximation~\cite[Corollary~4.8]{Jardine} with $A=\varnothing$ and $B=K$.  It produces
$K'$ by finitely many subdivision and dual-subdivision operations, a comparison $\gamma\colon K'\to K$, and
$\phi\colon K'\to X$ with $|\phi|\simeq g\circ|\gamma|$.  Thus $K'$ is finite and $v=|\gamma|$ is a homotopy
equivalence, as the cited theorem asserts for the subdivision comparisons.
\end{proof}

\begin{lemma}\label{lem:aperiodic}
If $L$ is aperiodic, then $M(Y,f)$ is aperiodic and both $BL$ and $BM(Y,f)$ are simply connected.
\end{lemma}
\begin{proof}
First consider $\theta\in O_{n+1}$.  For every $i\in[n]$, the sequence
$i,\theta(i),\theta^2(i),\ldots$ is monotone after its first step: it is nondecreasing if $i\leq\theta(i)$ and
nonincreasing if $\theta(i)\leq i$, because $\theta$ preserves order.  This sequence is therefore eventually constant.  Since
$[n]$ is finite, the maximum of the stabilization exponents for its vertices gives one $r$ with
$\theta^r=\theta^{r+1}$.  The same equality holds in the opposite monoid $O_{n+1}^{\op}$.

Now let $x=(y,l,i)$ belong to the ideal of $M(Y,f)$ and put $c=\hol_i(y)l$.  Induction using the first multiplication
formula in Definition~\ref{def:M} gives $x^r=(y,lc^{r-1},i)$ for all $r\geq1$.
Since $L$ is aperiodic, $c^s=c^{s+1}$ for some $s$, and consequently $x^{s+1}=x^{s+2}$.  Every element of $M(Y,f)$
therefore has stationary powers, so $M(Y,f)$ is aperiodic.

Finally, McDuff's lemma~\cite[Lemma~1]{McD} identifies $\pi_1(BL)$ with the universal group generated by the elements
of $L$, subject to the multiplication relations in $L$.  In a group, multiplying the relation $a^r=a^{r+1}$ by
$a^{-r}$ forces the image of $a$ to be the identity.  Since every element of $L$ has such a relation, every generator
of the universal group is trivial.  The same argument applies to the aperiodic monoid $M(Y,f)$.  Classifying spaces
of monoids are path connected, so both $BL$ and $BM(Y,f)$ are simply connected.
\end{proof}

\begin{lemma}\label{lem:reduced}
Suppose $L$ is aperiodic.  Let $Y'$ be a finite connected simplicial set and let $f'$ be a simplicial map
\[
 f'\colon Y'\longrightarrow N(L).
\]
Then there are a finite reduced simplicial set $Y$, a simplicial map $f\colon Y\to N(L)$ and a homotopy
equivalence $h\colon|Y'|\to|Y|$ with $|f|\circ h=|f'|$.
\end{lemma}
\begin{proof}
The reduction will be performed one edge at a time.  For each $a\in L$ we first construct a finite reduced
contractible simplicial set $D_a$, equipped with a map $\delta_a\colon D_a\to N(L)$ and a distinguished loop
$\varepsilon_a\colon\Delta^1\to D_a$ labelled by $a$.  Replacing an edge labelled by $a$ with this loop will identify
its two vertices without changing the homotopy type.

We use the following consequence of Lemma~\ref{lem:gluing} throughout.  If
$A\to X$ and $A\to B$ are maps of finite simplicial sets, one is a monomorphism, and $|A|\to|B|$ is a homotopy
equivalence, then $|X|\to|X\sqcup_A B|$ is a homotopy equivalence.  Indeed, realization preserves the pushout and
takes a simplicial monomorphism to a cofibration.

Choose $r$ with $a^r=a^{r+1}$ and put $e_a=a^r$.  Then $e_a^2=e_a$ and $e_a\,a\,e_a=e_a$.
Every power $a^s$ with $s\ge r$ is equal to $e_a$.  Notice that both $y=a$ and $y=1$ satisfy
$e_a\,y\,e_a=e_a$.  For any pair $e,y\in L$ satisfying $e^2=e$ and $e\,y\,e=e$, we now construct a finite contractible
simplicial set $H$ and a simplicial map $g_{e,y}\colon H\to N(L)$ with distinguished edges labelled by $e$,
$ye$, and $y$.
Let $H$ be obtained from an edge $\kappa\colon v_0\to v_1$ by two horn attachments.  The first adds a $2$-simplex
$\tau_1$ with $d_2\tau_1=d_1\tau_1=\kappa$ and $\mu:=d_0\tau_1$; the second adds a $2$-simplex $\tau_2$ with
$d_1\tau_2=\mu$, $d_0\tau_2=\kappa$, and $\beta:=d_2\tau_2$.
Here $\mu$ is a loop at $v_1$ and $\beta$ runs from $v_1$ to $v_0$.  Since $H$ is obtained from $\Delta^1$ by attaching
$\Delta^2$ along $\Lambda^2_0$ and then along $\Lambda^2_2$, the preceding observation applied twice shows that
$|H|\simeq|\Delta^1|$ and hence that $H$ is contractible.  The assignments
\[
 \kappa\longmapsto e,\qquad \mu\longmapsto ye,\qquad \beta\longmapsto y,\qquad
 \tau_1\longmapsto(e,ye),\qquad \tau_2\longmapsto(y,e)
\]
define a simplicial map $g_{e,y}\colon H\to N(L)$.  The two displayed pairs give the three prescribed edge
labels on each $2$-simplex; the only nonautomatic face identity is $e\,y\,e=e$.

For the fixed element $a$, let $H_a$ and $H_1$ be two copies of $H$, equipped with the maps $g_{e_a,a}$ and
$g_{e_a,1}$, respectively.  Use $H_1$ to form the simplicial quotient
$R_a:=H_1/\beta=H_1\mathop{\sqcup}_{\Delta^1}\Delta^0$.
Because $g_{e_a,1}(\beta)=1$, this map factors through $R_a\to N(L)$.  Collapsing the edge $\beta$, which joins
the two vertices, makes $R_a$ reduced.  The edge $\beta$ is a simplicial monomorphism, and
$|\Delta^1|\to|\Delta^0|$ is a homotopy equivalence; the observation gives $|H_1|\simeq|R_a|$, so $R_a$ is
contractible.  Write
$\bar\kappa\colon\Delta^1\to R_a$ for the image of $\kappa$, and set
$D_a:=H_a\mathop{\sqcup}_{\Delta^1}R_a$,
where the structure maps are $\kappa\colon\Delta^1\to H_a$ and
$\bar\kappa\colon\Delta^1\to R_a$.  The maps $g_{e_a,a}$ on $H_a$ and $R_a\to N(L)$ agree on the gluing edge,
whose label is $e_a$, and define $\delta_a\colon D_a\to N(L)$.  Gluing $\kappa$ to the loop $\bar\kappa$ identifies
$v_0$ and $v_1$, so $D_a$ is finite and reduced.  Since $\kappa$ is a monomorphism and
$|\bar\kappa|$ is a homotopy equivalence between contractible spaces, the observation gives
$|H_a|\simeq|D_a|$; thus $D_a$ is contractible.  The image
$\varepsilon_a$ of $\beta$ in $H_a$ is a loop with $\delta_a(\varepsilon_a)=a$.

We now reduce $Y'$ one vertex at a time.  The connected components of a simplicial set are already detected by its
$1$-skeleton.  Hence, if $Y'$ has more than one vertex, there is a nondegenerate edge $\gamma$ with distinct
vertices.  Put $a=f'(\gamma)$ and form the pushout
\[
 \begin{array}{ccc}
 \Delta^1&\xrightarrow{\ \varepsilon_a\ }&D_a\\
 {\scriptstyle \gamma}\big\downarrow&&\big\downarrow\\
 Y'&\longrightarrow&Y''.
 \end{array}
\]
The maps $f'$ and $\delta_a$ agree on $\Delta^1$, so they induce a simplicial map $f''\colon Y''\to N(L)$.  If
$j\colon Y'\to Y''$ denotes the pushout map, then $f''\circ j=f'$.  The map
$\gamma\colon\Delta^1\to Y'$ is a monomorphism: its edge is nondegenerate and its two vertices are distinct, so uniqueness
of degeneracy decompositions makes it injective in every simplicial degree.  On the other hand,
$|\varepsilon_a|$ is a homotopy equivalence because both spaces are contractible.  The observation therefore makes
$|j|\colon|Y'|\to|Y''|$ a homotopy equivalence.
Since $D_a$ has one vertex and the two vertices of $\gamma$ are identified, $Y''$ has exactly one fewer vertex than
$Y'$.  The simplicial set $Y''$ is connected because it is homotopy equivalent to $Y'$.  Iterating finitely many times produces a finite
reduced simplicial set $Y$, a simplicial map $f\colon Y\to N(L)$, and a simplicial map
$\varphi\colon Y'\to Y$ such that $f\circ\varphi=f'$ and $|\varphi|$ is a homotopy equivalence.  Take
$h=|\varphi|$.
\end{proof}

\begin{proof}[Proof of \ThmB]
\leavevmode\phantomsection\label{proof:main-B}%
Choose a finite connected simplicial set $K$ and a homotopy equivalence $u\colon|K|\to Y$
\cite[Theorem~2C.5]{Hatcher}.  Put $f_K=f\circ u\colon|K|\to BL$ and
$\widetilde f_K=\operatorname{rev}_L\circ f_K\colon|K|\to|N(L)|$.
By Lemma~\ref{lem:finite-simplicial-approximation}, there are a finite simplicial set $K_f$, a simplicial map
$\phi_f\colon K_f\to N(L)$, a homotopy equivalence $v_f\colon|K_f|\to|K|$, and a homotopy
$|\phi_f|\simeq\widetilde f_K\circ v_f$.  The simplicial set $K_f$ is connected because $v_f$ is a homotopy
equivalence and $K$ is connected.  Lemma~\ref{lem:reduced} supplies maps
$f_S\colon S(Y,f)\to N(L)$ and
$r_f\colon|K_f|\xrightarrow{\simeq}|S(Y,f)|$
such that $|f_S|\circ r_f=|\phi_f|$.  Choose a homotopy inverse
$\bar r_f\colon|S(Y,f)|\to|K_f|$ and define
\[
 h_{Y,f}:=u\circ v_f\circ\bar r_f\colon|S(Y,f)|\longrightarrow Y.
\]
After identifying $|N(L)|$ with $BL$ by $\operatorname{rev}_L^{-1}$, the preceding identities give
\[
 |f_S|\simeq |f_S|\circ r_f\circ\bar r_f
 =|\phi_f|\circ\bar r_f
 \simeq f_K\circ v_f\circ\bar r_f
 =f\circ h_{Y,f}.
\]
Set $M_f=M(S(Y,f),f_S)$.  Then \ThmA{} and
Lemma~\ref{lem:cofiber-invariance} give
\[
 BM_f\simeq\Cof(|f_S|)\simeq\Cof(f\circ h_{Y,f})\simeq\Cof(f).
\]
The monoid $M_f$ is finite by \ThmA{} and aperiodic by Lemma~\ref{lem:aperiodic}.
\end{proof}

\begin{lemma}\label{lem:finite-one-reduction}
Every finite simply connected CW complex is homotopy equivalent to a finite CW complex with one $0$-cell and no
$1$-cells.
\end{lemma}

\begin{proof}
Let $X$ be a finite simply connected CW complex.
Choose a maximal tree $T$ in the $1$-skeleton of $X$~\cite[Proposition~1A.1]{Hatcher}.  Since $T$ is a contractible
subcomplex, the quotient map $X\to X/T$ is a homotopy equivalence~\cite[Proposition~0.17]{Hatcher}.  Thus $X/T$ is a
finite simply connected CW complex with one $0$-cell.  The pair $(X/T,*)$ is $1$-connected, so the Cell Trading
Lemma~\cite[Proposition~4.2.1]{Geo}, with $n=1$ and $A=*$, replaces its finitely many $1$-cells by $3$-cells
relative to the vertex.  The resulting finite CW complex is homotopy equivalent to $X/T$ and has no $1$-cells.
\end{proof}

\begin{proof}[Proof of \ThmC]
\leavevmode\phantomsection\label{proof:main-C}%
Let $X$ be a finite simply connected CW complex.  By Lemma~\ref{lem:finite-one-reduction}, we may replace $X$ by a
homotopy-equivalent finite CW complex with one $0$-cell and no $1$-cells; continue to denote this model by $X$.
Attach its remaining cells one at a time,
in nondecreasing order of dimension, and write
\[
 \begin{gathered}
  *=X_0\subseteq X_1\subseteq\cdots\subseteq X_N=X,\\
  X_{j+1}=X_j\cup_{f_j}D^{m_j+1},
  \qquad f_j\colon S^{m_j}\longrightarrow X_j,
  \qquad m_j\geq1.
 \end{gathered}
\]
Since $X$ has no $1$-cells, $m_j\ge1$, so each attaching sphere is connected, as \ThmB{} requires.

Begin with the trivial aperiodic monoid $\Lambda_0=\{1\}$ and the canonical homotopy equivalence
$h_0\colon X_0=*\xrightarrow{\simeq}B\Lambda_0=*$. 
Suppose inductively that a finite aperiodic monoid $\Lambda_j$ and a homotopy
equivalence $h_j\colon X_j\to B\Lambda_j$ have been constructed.  Compose the next attaching map with $h_j$ to obtain
$\bar f_j=h_j\circ f_j\colon S^{m_j}\to B\Lambda_j$.
Because $S^{m_j}$ is a finite connected CW complex, \ThmB{} applies and
produces a finite aperiodic monoid $\Lambda_{j+1}$ together with a homotopy equivalence
$b_{j+1}\colon B\Lambda_{j+1}\xrightarrow{\simeq}\Cof(\bar f_j)$.
The definition of $\bar f_j$ gives the strictly commutative square
\[
\begin{array}{ccc}
 S^{m_j}&\xrightarrow{\ f_j\ }&X_j\\
 {\scriptstyle\id}\big\downarrow&&\big\downarrow{\scriptstyle h_j}\\
 S^{m_j}&\xrightarrow{\ \bar f_j\ }&B\Lambda_j.
\end{array}
\]
Both vertical maps are homotopy equivalences, so Lemma~\ref{lem:cofiber-invariance} gives a homotopy equivalence
$c_{j+1}\colon\Cof(f_j)\xrightarrow{\simeq}\Cof(\bar f_j)$.
Because $f_j$ is a map between CW complexes, its ordinary mapping cone is the cell attachment
$X_j\cup_{f_j}D^{m_j+1}=X_{j+1}$ and represents $\Cof(f_j)$.  The induction step is therefore the zigzag
\[
 X_{j+1}\xrightarrow{\ \simeq\ }\Cof(f_j)
 \xrightarrow[ c_{j+1}\ ]{\ \simeq\ }\Cof(\bar f_j)
 \xleftarrow[\ b_{j+1}\ ]{\ \simeq\ }B\Lambda_{j+1}.
\]
Choose $h_{j+1}\colon X_{j+1}\to B\Lambda_{j+1}$ in the homotopy class obtained by traversing this zigzag.  Because
$\Lambda_{j+1}$ is again finite and aperiodic, the step can be repeated for every cell.  At the end,
$B\Lambda_N\simeq X_N=X$, which proves the theorem.
\end{proof}

\begin{remark}
The monoids produced by this induction are very large, since each step uses all $n$-simplices of a subdivided sphere.
\ThmB{} is more flexible than the induction suggests: any decomposition of $X$ as an iterated mapping cone of maps
from finite connected complexes can be used.
\end{remark}

\section{Examples}\label{sec:cellular-examples}

The following examples make individual stages of the induction explicit.

\begin{example}\label{ex:cellular-spheres}
Let $S^q_\Delta=\Delta^q/\partial\Delta^q$, $q\geq1$.  This simplicial set has two $q$-simplices, the basepoint and
the top simplex, so $A_{q+1}:=M(S^q_\Delta,*)$ has
$\#A_{q+1}=\binom{2q+1}{q}+2(q+1)$ and $BA_{q+1}\simeq S^{q+1}$.
Write
$a_i=(*,i)$ and $b_i=(\sigma_q,i)$ for all $0\leq i\leq q$,
where $\sigma_q$ is the top simplex.  Then
$A_{q+1}=O_{q+1}^{\op}\sqcup\{a_0,\ldots,a_q,b_0,\ldots,b_q\}$, and, for $u,v\in\{a,b\}$,
\[
 u_i v_j=u_j,
 \qquad
 u_i\theta^{\op}=u_{\theta(i)},
 \qquad
 \theta^{\op}a_i=a_i,
 \qquad
 \theta^{\op}b_i=
 \begin{cases}
  b_i,&\theta=\id,\\
  a_i,&\theta\ne\id.
 \end{cases}
\]
The ideal on the $a_i,b_i$ is a $2\times(q+1)$ rectangular band.  The orders of
$A_2,A_3,A_4$ are respectively $7,16,43$, and their classifying spaces are respectively
$S^2,S^3,S^4$ up to homotopy.

More generally, for $Y=\bigvee^r S^q_\Delta$, one has
\[
 \#M(Y,*)=\binom{2q+1}{q}+(q+1)(r+1),
 \qquad BM(Y,*)\simeq\bigvee^r S^{q+1}.
\]
For $q=1$ and $r=2$, this gives a nine-element aperiodic monoid whose classifying space is
$S^2\vee S^2$.
\end{example}

\begin{example}
Let $C_2$ be the cyclic group of order two, let $Y=\sk_2N(C_2)$ be the simplicial $2$-skeleton, and use
$L=\{1\}$.  Then $|Y|\cong\RP^2$, $\#Y_2=4$, $\#M(Y,*)=10+3\cdot4=22$, and
$BM(Y,*)\simeq\Sigma\RP^2\simeq S^2\cup_2 e^3$.
\end{example}

\begin{example}\label{ex:filling-sphere}
This is the simply connected analogue of Example~\ref{ex:holonomy}, and a warm-up for
Example~\ref{ex:cp2}.  Let $L_0$ be the five-element monoid of Example~\ref{ex:intro-sphere}, so
that $BL_0\simeq S^2$, and let $u\colon L_0\times L_0\to\mathbb Z$ be the normalized $2$-cocycle
\[
 u(x_{ij},x_{k\ell})=
 \begin{cases}
  1,&j=k=0,\\
  0,&\text{otherwise},
 \end{cases}
 \qquad
 u(1,a)=u(a,1)=0.
\]
In the normalized bar complex, the chain
\[
 z=[x_{00}\mid x_{00}]-[x_{00}\mid x_{10}]-[x_{01}\mid x_{00}]+[x_{01}\mid x_{10}]
\]
is a $2$-cycle with $u(z)=1$.  Since $H_1(BL_0;\mathbb Z)=0$ and $H_2(BL_0;\mathbb Z)\cong\mathbb Z$
by Example~\ref{ex:intro-sphere}, $[u]$ generates
$H^2(BL_0;\mathbb Z)\cong\operatorname{Hom}(H_2(BL_0;\mathbb Z),\mathbb Z)$ and $[z]$ generates
$H_2(BL_0;\mathbb Z)$.

Let $Y_z\subseteq N(L_0)$ be the simplicial subset generated by the four $2$-simplices
\[
 \sigma_1=(x_{00},x_{00}),\quad \sigma_2=(x_{00},x_{10}),\quad
 \sigma_3=(x_{01},x_{00}),\quad \sigma_4=(x_{01},x_{10})
\]
occurring in $z$, and let $f_z\colon Y_z\hookrightarrow N(L_0)$ be the inclusion.  The faces
$d_0(p,q)=q$, $d_1(p,q)=pq$ and $d_2(p,q)=p$ show that $Y_z$ is reduced, with nondegenerate edges
$x_{00},x_{10},x_{01}$ and nondegenerate $2$-simplices $\sigma_1,\dots,\sigma_4$.  In the normalized
chain complex,
\[
 \partial\sigma_1=x_{00},\qquad \partial\sigma_2=x_{10},\qquad
 \partial\sigma_3=x_{01},\qquad \partial\sigma_4=x_{10}-x_{00}+x_{01}.
\]
These formulas give $\ker\partial_2=\mathbb Z\cdot z$; since $Y_z$ has no nondegenerate simplices
above dimension $2$, it follows that $H_2(Y_z;\mathbb Z)=\mathbb Z\cdot z$.  The group relations
$[p][q]=[pq]$ associated to $\sigma_1,\sigma_2,\sigma_3$ kill the three edge generators, so
$\pi_1(|Y_z|)=1$ and $H_1(Y_z;\mathbb Z)=0$.  The Hurewicz isomorphism supplies a map
$S^2\to|Y_z|$ representing $z$; it is a homology equivalence and hence a homotopy equivalence by
the homological Whitehead theorem.  Thus $|Y_z|$ is a homotopy $2$-sphere with fundamental class
$z$.  With the orientations determined by $z$ and $[u]$, the equality $u(z)=1$ shows that
$|f_z|\colon|Y_z|\to BL_0$ has degree one.

Here $n=2$, and $(Y_z)_2$ has eleven elements, four nondegenerate and seven degenerate, so
$\#M(Y_z,f_z)=\binom52+3\cdot5\cdot11=175$.  For $\sigma=(p,q)\in(Y_z)_2$ the holonomies are
$\hol_0(\sigma)=pq$, $\hol_1(\sigma)=q$ and $\hol_2(\sigma)=1$.  As in
Example~\ref{ex:holonomy}, the constant maps $c_k\in O_3$ act by
\[
 c_k^{\op}\,(\sigma,1,i)=(*_2,\hol_k(\sigma),i)=\lambda_i(\hol_k(\sigma)),\qquad k=0,1,2.
\]
The vertex operators send $(\sigma,1,i)$ to $\lambda_i(pq)$, $\lambda_i(q)$ and
$\lambda_i(1)=e_i$.  These are the three vertices of $(\sigma,1)\in(E_{f_z})_2$, read in the
corner $e_iM(Y_z,f_z)e_i\cong L_0$.  By \ThmA, $BM(Y_z,f_z)\simeq\Cof(|f_z|)$ is contractible,
since $|f_z|$ is a degree-one map between homotopy $2$-spheres.

Now replace $f_z$ by the constant map $*$, so that every edge of $Y_z$ has label $1$.  Every
holonomy becomes $1$ and every $c_k^{\op}(\sigma,1,i)$ equals $e_i$, and \ThmA{} gives
$BM(Y_z,*)\simeq BL_0\vee\Sigma|Y_z|\simeq S^2\vee S^3$.  The monoids $M(Y_z,f_z)$ and
$M(Y_z,*)$ are aperiodic by Lemma~\ref{lem:aperiodic}, have the same $175$ elements, and differ
only in their holonomies.  That difference decides whether the new $3$-cell fills the $2$-sphere
$BL_0$ or is attached along a constant map.
\end{example}

\begin{example}\label{ex:cp2}
Let $L_0$ be the five-element aperiodic monoid of
Example~\ref{ex:intro-sphere}, so $BL_0\simeq S^2$.  Recall that the standard CW structure on $\CP^2$ is obtained
from $S^2$ by attaching one $4$-cell along the Hopf map $\eta\colon S^3\to S^2$; thus
$\CP^2\simeq S^2\cup_\eta D^4=\Cof(\eta)$.
We will realize this cell attachment using \ThmA.  The monoid $L_0$ models the $2$-skeleton $S^2$; we next construct
a finite reduced simplicial set $Y_\eta$ modeling the attaching sphere $S^3$ and a simplicial map
$f_\eta\colon Y_\eta\to N(L_0)$ whose realization represents $\eta$.  Let $u$ and $z$ be as in
Example~\ref{ex:filling-sphere}, and fix the equivalence $BL_0\simeq S^2$ so that $[u]$ corresponds to the standard
positive generator of $H^2(S^2;\mathbb Z)$.
Form the central extension $\widehat L_0=L_0\times_u\mathbb Z$, with multiplication
\[
 (a,r)(b,s)=(ab,r+s+u(a,b)).
\]
The computational supplement is available in the
\href{https://github.com/OrenXGitHub/monoid-computations}{GitHub repository}.  The program
\href{https://github.com/OrenXGitHub/monoid-computations/blob/main/cofiber_monoids.py}{\texttt{cofiber\_\allowbreak monoids.py}}
contains both the finite input and the code that verifies it; the repository as a whole is described in the
\hyperref[app:computational-supplement]{Appendix}.  The input consists of twelve explicitly listed triples
$(p_1,p_2,p_3)\in\widehat L_0^3=N(\widehat L_0)_3$.  Regarded as nondegenerate $3$-simplices, they generate a
finite reduced simplicial subset $Y_\eta\subseteq N(\widehat L_0)$.  The program closes them under faces and
degeneracies, obtaining one vertex, seven
nondegenerate edges, eighteen nondegenerate $2$-simplices, twelve nondegenerate $3$-simplices, and
$\#(Y_\eta)_3=88$ simplices in degree $3$ when degeneracies are included.  Projection $\widehat L_0\to L_0$ induces the
simplicial map $f_\eta\colon Y_\eta\to N(L_0)$.

The exact calculations described in the \hyperref[app:computational-supplement]{Appendix} prove
\[
 \pi_1(|Y_\eta|)=1,\qquad
 H_i(Y_\eta;\mathbb Z)\cong
 \begin{cases}
  \mathbb Z,&i=0,3,\\
  0,&i\ne0,3,
 \end{cases}
 \qquad H(|f_\eta|)=1.
\]
The first two assertions make $|Y_\eta|$ a homotopy $3$-sphere: Hurewicz gives $\pi_2=0$ and then
$\pi_3\cong H_3\cong\mathbb Z$, so a generator $S^3\to|Y_\eta|$ is a homology equivalence and hence a homotopy
equivalence by the homological Whitehead theorem.  With the orientation fixed above, the last assertion says that
$|f_\eta|$ represents the Hopf map $\eta\in\pi_3(S^2)$.  The program also constructs the complete multiplication
table of $M(Y_\eta,f_\eta)$ and checks its identity and associativity laws; its finite-range cohomology calculation is
an independent consistency check.

Define $M_{\CP^2}:=M(Y_\eta,f_\eta)$.
Since $n=3$, \ThmA{} gives
\[
 \begin{aligned}
 \#M_{\CP^2}
 &=\binom{7}{3}+4(\#L_0)\#\bigl((Y_\eta)_3\bigr)
 =35+4\cdot5\cdot88=1795,\\
 BM_{\CP^2}
 &\simeq\Cof(|f_\eta|)\simeq\Cof(\eta)\simeq\CP^2.
 \end{aligned}
\]
$M_{\CP^2}$ is an explicit finite instance of the induction step; it is not intended to minimize the order of the
monoid.  A separate $122$-element monoid $M_{122}$ with $BM_{122}\simeq\CP^2$ is described in the
\hyperref[app:computational-supplement]{Appendix}.
Let $x$ generate the one-dimensional vector space $H^2(\CP^2;\mathbb F_2)$.
Then
\[
 x^2=Sq^2(x)\ne0,
\]
which distinguishes this classifying space from $S^2\vee S^4$.  This has a direct
monoid-algebraic meaning: for $M=M_{\CP^2}$ and $k=\mathbb F_2$, the bar resolution identifies
$H^*(BM;k)$ with $\operatorname{Ext}_{k[M]}^*(k,k)$ as rings.  If $\omega\colon M^2\to k$ represents $x$, then
$x^2$ is represented by
\[
 (a,b,c,d)\longmapsto\omega(a,b)\omega(c,d).
\]
This cochain represents the Yoneda square of the corresponding degree-$2$ extension, so that square is nonzero.  In this sense the multiplication table detects the Hopf attachment.
\end{example}

\section{The finite fundamental-group question}\label{sec:fundamental-group}

All monoids constructed in the preceding sections are aperiodic.  We now drop that condition, but even arbitrary finite monoids cannot
realize every connected finite CW homotopy type: $\pi_1(BM)$ is the universal group $U(M)$ defined in the
introduction, and $U(M)$ is finite when $M$ is finite.  Let $\eta_M\colon M\to U(M)$ be the universal homomorphism.
For each $a\in M$ there are $r<s$ with $a^r=a^s$, and cancellation in $U(M)$ gives
$\eta_M(a)^{s-r}=1$.  Hence $\eta_M(a)^{-1}=\eta_M(a^{s-r-1})$, with exponent zero interpreted as the identity.
Replace every inverse in a group word by this nonnegative power.  The resulting positive word is $\eta_M(b)$ for
the product $b\in M$ of its letters.  Thus $\eta_M$ is surjective and $U(M)$ is
finite.

This necessary condition suggests the following conjecture.

\begin{conjecture}\label{conj:finite-fundamental-group}
Let $X$ be a connected finite CW complex.  If $\pi_1(X)$ is finite, then
$X\simeq BM$ for some finite monoid $M$.
\end{conjecture}

Together with the necessary condition above, Conjecture~\ref{conj:finite-fundamental-group} would say that the
connected finite CW homotopy types realized by finite monoids are exactly those with finite fundamental group.
\ThmC{} proves the simply connected case.  The monoid in
Conjecture~\ref{conj:finite-fundamental-group} cannot be aperiodic unless $X$ is simply connected, by
Lemma~\ref{lem:aperiodic}.

\begin{remark}
There is a concise poset description of every pointed connected finite CW homotopy type, but it retains the
fundamental group separately and therefore does not produce the finite monoid sought in
Conjecture~\ref{conj:finite-fundamental-group}.  Choose a finite simplicial complex $K$ homotopy equivalent to $X$,
let $G=\pi_1(|K|)$, and let $P$ be the face poset of the universal cover $\widetilde K$.  The deck action of $G$ on
$P$ is free, $|N(P)|\cong|\Sd\widetilde K|$ is simply connected, and
$N(P)/G\cong\Sd K$ is finite.  Hence
$|N(P)|/G\cong|\Sd K|\simeq X$, and $|N(P)|\to|N(P)|/G$ is the universal covering with deck group $G$.
When $G$ is finite, $P$ is finite as well, so $X$ is described by a finite poset with a free action of a finite
group; Conjecture~\ref{conj:finite-fundamental-group} asks whether this description can always be replaced by a
single finite monoid.

Conversely, suppose a discrete group $G$ acts freely on a poset $P$, $|N(P)|$ is simply connected, and $N(P)/G$
has finitely many nondegenerate simplices.  For each vertex of $N(P)$, only finitely many translates of the finitely
many orbit representatives of simplices can meet its star; thus $N(P)$ is locally finite.  No nonidentity element
stabilizes a chain, because an order-preserving permutation of a finite chain fixes every vertex.  The action is
therefore a covering action in the sense of~\cite[Proposition~1.40]{Hatcher}, so the quotient is a connected finite
CW complex with universal cover $|N(P)|$ and fundamental group $G$.  If
$A=\pi_2(|N(P)|)$, the action gives the monodromy $G$-module structure on $A$, and the remaining datum in the
$2$-type is the Postnikov class in $H^3(G;A)$.
\end{remark}

\medskip
\noindent\textbf{Acknowledgements.} I would like to thank Clark Barwick and David Blanc for helpful conversations about this work.

\medskip
\noindent\textbf{AI declaration.}\par
The author takes responsibility for the contents of this article and the
GitHub repository cited below.  The key idea of realizing $\mathbb{CP}^2$ by a finite monoid arose in a conversation between the author
and Anthropic's Claude Opus 5.5.  Starting from that idea, the author supplied
all of the theoretical arguments in the article.  All of the content in the GitHub
repository was written by Claude.

\appendix
\phantomsection
\section*{Appendix}\label{app:computational-supplement}
\addcontentsline{toc}{section}{Appendix}

The repository available at
\begin{center}
\url{https://github.com/OrenXGitHub/monoid-computations}
\end{center}
contains the Python source, documentation, and dependency information for the examples.  Its file
\href{https://github.com/OrenXGitHub/monoid-computations/blob/main/cofiber_monoids.py}{\texttt{cofiber\_monoids.py}}
implements $M(Y,f)$, records the twelve generators of $Y_\eta$, and performs the following exact calculations.
For the normalized chain complex, let
$D_i\colon C_i(Y_\eta;\mathbb Z)\to C_{i-1}(Y_\eta;\mathbb Z)$ be the boundary.  The program verifies
$D_2D_3=0$ and exhibits determinant $-1$ minors of sizes $7\times7$ in
$D_2\colon\mathbb Z^{18}\to\mathbb Z^7$ and $11\times11$ in
$D_3\colon\mathbb Z^{12}\to\mathbb Z^{18}$.  Thus $D_2$ has rank $7$, while $D_3$ has rank $11$ and primitive
image.  It follows that $H_1=0$, that $\ker D_2$ and $\operatorname{im}D_3$ have rank $11$, and that
$H_2=\ker D_2/\operatorname{im}D_3$ is a finite subgroup of the torsion-free group
$\mathbb Z^{18}/\operatorname{im}D_3$; hence $H_2=0$.  Finally $H_3=\ker D_3\cong\mathbb Z$.
Using one generator for each of the seven nondegenerate edges and the triangle relation
$[d_2\tau][d_0\tau]=[d_1\tau]$ for every nondegenerate $2$-simplex, the program reduces all generators to the
identity and thereby proves $\pi_1(|Y_\eta|)=1$.

For the Hopf invariant, write $\pi\colon\widehat L_0\to L_0$ for projection and
$\nu(a,r)=r$.  The extension law gives $d(-\nu)=\pi^*u$.  The program solves exactly for a normalized integral
$3$-cochain $w$ with $dw=u\smile u$ and forms
$\Theta=(-\nu)\smile\pi^*u-\pi^*w$, so that
\[
 \Theta((a_1,r_1),(a_2,r_2),(a_3,r_3))
 =-r_1u(a_2,a_3)-w(a_1,a_2,a_3).
\]
Its restriction to $Y_\eta$ is a cocycle, and its value on the computed fundamental cycle is $1$.  By the cochain
definition of the Hopf invariant, this proves $H(|f_\eta|)=1$.  The same file constructs the full
$1795$-element multiplication table, verifies the monoid laws, and computes its mod-$2$ cohomology ring through
degree $9$ as an independent check.

The file
\href{https://github.com/OrenXGitHub/monoid-computations/blob/main/annulus_cp2.py}{\texttt{annulus\_cp2.py}}
writes the multiplication table of a separate $122$-element monoid $M_{122}$.  The stand-alone note
\href{https://github.com/OrenXGitHub/monoid-computations/blob/main/ANNULUS_CP2_PROOF.md}{\texttt{ANNULUS\_CP2\_\allowbreak PROOF.md}}
proves $BM_{122}\simeq\CP^2$ by identifying the relevant category-of-elements map with the Hopf fibration.  Its
finite categorical, simplicial, and integral-cocycle calculations are checked by
\href{https://github.com/OrenXGitHub/monoid-computations/blob/main/verify_annulus_cp2.py}{\texttt{verify\_annulus\_\allowbreak cp2.py}}.
The program \texttt{annulus\_cp2.py} separately checks the monoid laws and computes
$H^*(BM_{122};\mathbb F_p)$ for $p=2,3,5,7$.  The file
\href{https://github.com/OrenXGitHub/monoid-computations/blob/main/check_cp2_corners.py}{\texttt{check\_cp2\_\allowbreak corners.py}}
tests proper idempotent corners of both models.  The \texttt{README.md} gives commands and expected output, and
\texttt{requirements.txt} records the Python~3.10 and NumPy~1.24 requirements.  We do not know whether $122$ is
minimal; determining the least order of a finite monoid $M$ with $BM\simeq\CP^2$ remains open.  None of the general
theorems depends on these computations.


\begin{thebibliography}{10}
\bibitem{Burke} N.~Burke, \emph{Homotopy theory of monoids and group completion}, Ph.D. thesis, University of
Cambridge, 2021, \href{https://www.repository.cam.ac.uk/handle/1810/325358}{repository}.
\bibitem{Fie84} Z.~Fiedorowicz, \emph{Classifying spaces of topological monoids and categories}, Amer. J. Math. \textbf{106}
(1984), 301--350, \href{https://doi.org/10.2307/2374307}{DOI}.
\bibitem{Fie02} Z.~Fiedorowicz, \emph{A counterexample to a group completion conjecture of J.~C.~Moore}, Algebr. Geom.
Topol. \textbf{2} (2002), 33--35, \href{https://doi.org/10.2140/agt.2002.2.33}{DOI}.
\bibitem{Geo} R.~Geoghegan, \emph{Topological methods in group theory}, Graduate Texts in Mathematics, vol.~243,
Springer, New York, 2008, \href{https://doi.org/10.1007/978-0-387-74614-2}{DOI}.
\bibitem{GZ} P.~Gabriel and M.~Zisman, \emph{Calculus of fractions and homotopy theory}, Springer, 1967,
\href{https://doi.org/10.1007/978-3-642-85844-4}{DOI}.
\bibitem{Hatcher} A.~Hatcher, \emph{Algebraic topology}, Cambridge University Press, 2002,
\href{https://pi.math.cornell.edu/~hatcher/AT/AT.pdf}{online edition}.
\bibitem{Jardine} J.~F.~Jardine, \emph{Simplicial approximation}, Theory Appl. Categ. \textbf{12} (2004), 34--72,
\href{https://www.tac.mta.ca/tac/volumes/12/2/12-02.pdf}{online}.
\bibitem{May} J.~P.~May, \emph{A concise course in algebraic topology}, University of Chicago Press, Chicago, 1999,
\href{https://www.math.uchicago.edu/~may/CONCISE/ConciseRevised.pdf}{online edition}.
\bibitem{MS} A.~Martinian and B.~Steinberg, \emph{Realizing wedges of Moore spaces as classifying spaces of finite
semigroups}, Semigroup Forum (2026), \href{https://doi.org/10.1007/s00233-026-10648-0}{DOI}.
\bibitem{McD} D.~McDuff, \emph{On the classifying spaces of discrete monoids}, Topology \textbf{18} (1979), 313--320,
\href{https://doi.org/10.1016/0040-9383(79)90022-3}{DOI}.
\bibitem{Qui} D.~Quillen, \emph{Higher algebraic K-theory.~I}, Lecture Notes in Math. \textbf{341}, Springer, 1973,
85--147, \href{https://doi.org/10.1007/BFb0067053}{DOI}.
\bibitem{RS} J.~Rhodes and B.~Steinberg, \emph{The $q$-theory of finite semigroups}, Springer Monographs in
Mathematics, Springer, New York, 2009.
\bibitem{Swe} D.~Sweeney, \emph{Homotopy types of small semigroups}, preprint, 2025,
\href{https://arxiv.org/abs/2502.05276}{arXiv:2502.05276}.
\bibitem{Tho} R.~W.~Thomason, \emph{Homotopy colimits in the category of small categories}, Math. Proc. Cambridge Philos.
Soc. \textbf{85} (1979), 91--109,
\href{https://doi.org/10.1017/S0305004100055535}{DOI}.
\end{thebibliography}
\end{document}